\documentclass[a4paper, 12pt, tikz]{emilsart}
\title{Permutation patterns and higher connectivity in inversion graphs}
\author{Emil Verkama}
\affiliationblock{Emil Verkama}{Department of Mathematics, KTH Royal Institute of Technology, Stockholm, Sweden}{verkama@kth.se}

\DeclareMathOperator{\Sub}{Sub}
\DeclareMathOperator{\deficit}{deficit}

\usepackage[backend=biber, doi=false, isbn=false, giveninits=true, style=alphabetic, maxbibnames=10]{biblatex}
\DeclareFieldFormat*{title}{``#1''}
\begin{document}
\maketitle

\begin{abstract}
  We prove that any permutation \(\pi \in S_n\) whose inversion graph is \(k\)-connected has at least \(\frac{k+1}{2}n - C\) inversions, where \(C\) depends only on \(k\). In particular, inversion graphs beat the general lower bound \(\frac{kn}{2}\) for the number of edges in a \(k\)-connected graph. Furthermore, we characterize the permutation patterns \(\sigma\) having the property that if \(\pi\) avoids \(\sigma\), then the bound improves to \(kn - C\). The collection of patterns is \(\Av(321, 2341, 4123)\), and it is independent of \(k \geq 2\).
  
  This work is motivated by a result of Linusson and Verkama (2025) that can be reinterpreted as stating that any \(1324\)-avoiding permutation whose inversion graph satisfies a certain weaker notion of \(2\)-connectivity must have at least \(2n - C\) inversions. Our final result gives a characterization of the patterns having this stronger property.
\end{abstract}

\begingroup
\hypersetup{linkcolor=black}
\tableofcontents
\endgroup

\section{Introduction}

The inversion graph \(G_{\pi}\) (also called the permutation graph) of a permutation \(\pi \in S_n\) has vertices labeled by the entries of \(\pi\) and edges between any two entries forming an inversion: the one on the left is larger than the one on the right. This is a classical, well-studied class of graphs. Inversion graphs are, for example, characterized by the property that both \(G_{\pi}\) and \(\overline{G_{\pi}}\) are comparability graphs of posets \cite{pnueli_transitive_1971}. See also \cite{even_permutation_1972} and \cite{golumbic_algorithmic_2004}.

Our perspective comes from enumerative combinatorics: several properties of \(\pi\) that are of enumerative interest translate nicely to \(G_{\pi}\). For instance, the components of the direct-sum decomposition of \(\pi\), which are of high importance in the study of pattern-avoiding permutations, correspond exactly to the connected components of \(G_{\pi}\). Further, patterns in \(\pi\) become (labeled) induced subgraphs of \(G_{\pi}\). 

In particular, a permutation is indecomposable if and only if its inversion graph is connected. The study of indecomposability and permutation patterns using inversion graphs is not new; see e.g.\ \cite{acan_random_2016,acan_connected_2013,koh_connected_2007} and \cite{alecu_letter_2022,atminas_well-quasi-order_2015,bassino_linear-sized_2022,brignall_labelled_2022}.

However, \emph{higher} connectivity of \(G_{\pi}\) does not appear to have been studied much. Our first result is that inversion graphs beat the general bound \(|E(G)| \geq \frac{kn}{2}\) for \(k\)-connected graphs \(G\). The general bound is known to be tight \cite{harary_maximum_1962} for every \(k \geq 2\).

\begin{theorema} \label{thm:k connected}
  Fix \(k \geq 1\). There exists a constant \(C\) depending only on \(k\) such that every permutation \(\pi\) whose inversion graph is \(k\)-connected satisfies
  \begin{equation*}
    \inv(\pi) \geq \frac{k+1}{2} |\pi| - C.
  \end{equation*}
\end{theorema}

The case \(k = 2\) (with the explicit constant \(C = 2\) from Section~\ref{sec:2-connected}) proves a conjecture of Verkama \cite[Conjecture~5.0.4]{verkama_inversions_2026}. The crucial tool is that inversion graphs are \(4\)-chordal: they contain no induced cycles of length five or longer \cite[Proposition~2.3]{tenner_cyclomatic_2026}. Our Proposition~\ref{prop:l chordal bound} in Section~\ref{sec:2-connected} generalizes a corresponding edge bound for bipartite \(4\)-chordal graphs of Zhang and Huang \cite{zhang_minimum_2024}.

The initial motivation behind this work is a result of Linusson and Verkama \cite[Theorem~8]{linusson_enumerating_2025} that every \(1324\)-avoiding permutation \(\pi \in S_n\) (this is a particularly important class) that has at most \(2n - 7\) inversions must be decomposable or `almost' decomposable. That result was useful in making progress on a conjecture posed by Claesson, Jelínek and Steingrímsson \cite{claesson_upper_2012}, as decomposable \(1324\)-avoiders have a particularly simple structure. However, the proof of \cite[Theorem~8]{linusson_enumerating_2025} is obscure. Later, Verkama generalized the result to the patterns \(1243\), \(1342\) and \(2143\) \cite[Note~C]{verkama_inversions_2026}, using mostly the same strategy. Claesson, Linusson, Ulfarsson and Verkama \cite{claesson_inversion_2026} and Meng \cite{meng_residual_2026} have expanded on the result in other ways.

We observe that a permutation \(\pi \in S_n\) is neither decomposable nor almost decomposable if and only if its inversion graph \(G_{\pi}\) satisfies a certain weaker notion of \(2\)-connectivity: each of the graphs \(G_{\pi}\), \(G_{\pi} \setminus 1\), \(G_{\pi} \setminus n\), \(G_{\pi} \setminus \pi(1)\) and \(G_{\pi} \setminus \pi(n)\) is connected. The result of Linusson and Verkama is that each such \emph{weakly \(2\)-connected} \(1324\)-avoiding inversion graph must have at least \(2n - 6\) edges.

This is a reason to study lower bounds on the number of edges in \(G_{\pi}\) under assumptions on connectivity and pattern avoidance. Theorem~\ref{thm:k connected} above shows that under normal \(k\)-connectivity, we do not even need pattern avoidance to get a nontrivial bound. The following theorem, which is our main result, characterizes exactly which patterns \(\sigma\) have a property analogous to that of \(1324\) in the work of Linusson and Verkama: if \(\pi \in S_n\) has a \(k\)-connected inversion graph and it avoids \(\sigma\), then \(\pi\) must have at least \(kn - C\) inversions for some constant \(C\). Write \(\Av(\sigma)\) for the set of all permutations avoiding \(\sigma\).

\begin{theorema} \label{thm:k connected av}
  Fix \(k \geq 2\), and let \(\sigma\) be a pattern. The following are equivalent:
  \begin{enumerate}[(i)]
    \item \label{thm:k connected av:forcing} There exists a constant \(C\) depending only on \(k\) and \(\sigma\) such that every permutation \(\pi \in \Av(\sigma)\) whose inversion graph is \(k\)-connected satisfies
    \begin{equation*}
      \inv(\pi) \geq k|\pi| - C.
    \end{equation*}
    \item \label{thm:k connected av:patterns} We have \(\sigma \in \Av(321, 2341, 4123)\).
    \item \label{thm:k connected av:graphs} The components of \(G_{\sigma}\) are paths and \(4\)-cycles.
  \end{enumerate}
\end{theorema}

Notice the minor miracle that conditions \ref{thm:k connected av:patterns} and \ref{thm:k connected av:graphs} are independent of \(k\).

This result already sheds some light on \cite[Theorem~8]{linusson_enumerating_2025} and \cite[Note~C]{verkama_inversions_2026}: each of the patterns \(1324\), \(1243\), \(1342\) and \(2143\) is, of course, in the class \(\Av(321, 2341, 4123)\). This theorem is not, however, a generalization of the aforementioned results, since we assume normal \(2\)-connectivity (setting \(k = 2\)) instead of the weak form. 

Our last task is, thus, to apply our new perspective and find an honest generalization of \cite[Theorem~8]{linusson_enumerating_2025}. In fact, as in the result above, we are able to get a full characterization. Some notation: let the group \(\Gamma = \mathbb Z_2 \times \mathbb Z_2\) act on the set of all permutations by the inverse \(\pi \mapsto \pi^{-1}\) and the reverse-complement \(\pi \mapsto \pi^{\rc}\); these maps induce isomorphisms of the inversion graphs. The orbit of a set \(B\) of permutations under this action will be denoted by \(\Gamma B\). Furthermore, denote by \(\Sub(B)\) the set of all patterns contained in the elements of \(B\), and let \(\id_t\) be the identity permutation of size \(t\).

\begin{theorema} \label{thm:weakly 2-connected}
  Let \(\sigma\) be a pattern. The following are equivalent:
  \begin{enumerate}[(i)]
    \item\label{thm:weakly 2-connected:forcing} There exists a constant \(C\) depending only on \(\sigma\) such that every \(\pi \in \Av(\sigma)\) whose inversion graph is weakly \(2\)-connected satisfies
    \begin{equation*}
      \inv(\pi) \geq 2|\pi| - C.
    \end{equation*}
    \item\label{thm:weakly 2-connected:patterns} We have \(\sigma \in \Sub(\Gamma B)\), where
    \begin{equation*}
      B = \left\{1324, 1342, 2413\right\} \cup \left\{21 \oplus {\id_t} \oplus 21 : t \in \mathbb Z_{\geq 0}\right\}.
    \end{equation*}
  \end{enumerate}
\end{theorema}

\paragraph{Outline}

We begin by building the theory of \(2\)-connected inversion graphs and their ear decompositions in Section~\ref{sec:2-connected}. This leads to a simple proof (and generalization) of the \(k = 2\) case of Theorem~\ref{thm:k connected}. In Section~\ref{sec:2-connected av}, we continue with some deeper structural observations regarding specialized ear decompositions of inversion graphs. The resulting structural understanding is precise enough to handle pattern avoidance, and the \(k = 2\) case of Theorem~\ref{thm:k connected av} follows.

In Section~\ref{sec:menger}, we prove the full versions of both Theorem~\ref{thm:k connected} and Theorem~\ref{thm:k connected av}. The surprise is that these statements contain very little additional depth over the \(k = 2\) versions we already proved. The only new ingredient is a particular application of Menger's theorem. 

Section~\ref{sec:weak 2-connectivity} concerns Theorem~\ref{thm:weakly 2-connected}. We develop a strategy that is reminiscent of that for Theorem~\ref{thm:k connected av}, adapted to weak \(2\)-connectivity. Finally, Section~\ref{sec:future work} proposes possible directions for future research, and Section~\ref{sec:appendix} in the appendix collects certain useful constructions.

\paragraph{Preliminaries} 

\begin{figure}[t]
  \centering
  \begin{tikzpicture}[scale=0.6]
    \foreach \x/\y in {1/6, 2/2, 3/4, 4/1, 5/7, 6/8, 7/3, 8/5}
      \coordinate (v\y) at (\x,\y);
    \foreach \u/\v in {6/2, 6/4, 6/1, 6/3, 6/5, 2/1, 4/1, 4/3, 7/3, 7/5, 8/3, 8/5}
      \draw[rpath] (v\u) -- (v\v);
    \foreach \y in {1, ..., 8}
      \node[dot] at (v\y) {};
    \foreach \y/\p in {1/below right, 2/below left, 3/below right, 4/above, 5/below right, 6/above left, 7/above left, 8/above}
      \node[\p] at (v\y) {\(\y\)};

    \begin{scope}[shift={(10,0)}]
      \foreach \x/\y in {1/2, 2/3, 3/1, 4/4, 5/7, 6/5, 7/8, 8/6}
        \coordinate (v\y) at (\x,\y);
      \foreach \u/\v in {2/1, 3/1, 7/5, 7/6, 8/6}
        \draw[rpath] (v\u) -- (v\v);
      \foreach \y in {1, ..., 8}
        \node[dot] at (v\y) {};
      \foreach \y/\p in {1/below, 2/above left, 3/above, 4/above, 5/below, 6/right, 7/above left, 8/above}
        \node[\p] at (v\y) {\(\y\)};
    \end{scope}
  \end{tikzpicture}
  \caption{The inversion graphs of \(62417835\) (left) and \(23147586\) (right).}
  \label{fig:inv graph example}
\end{figure}
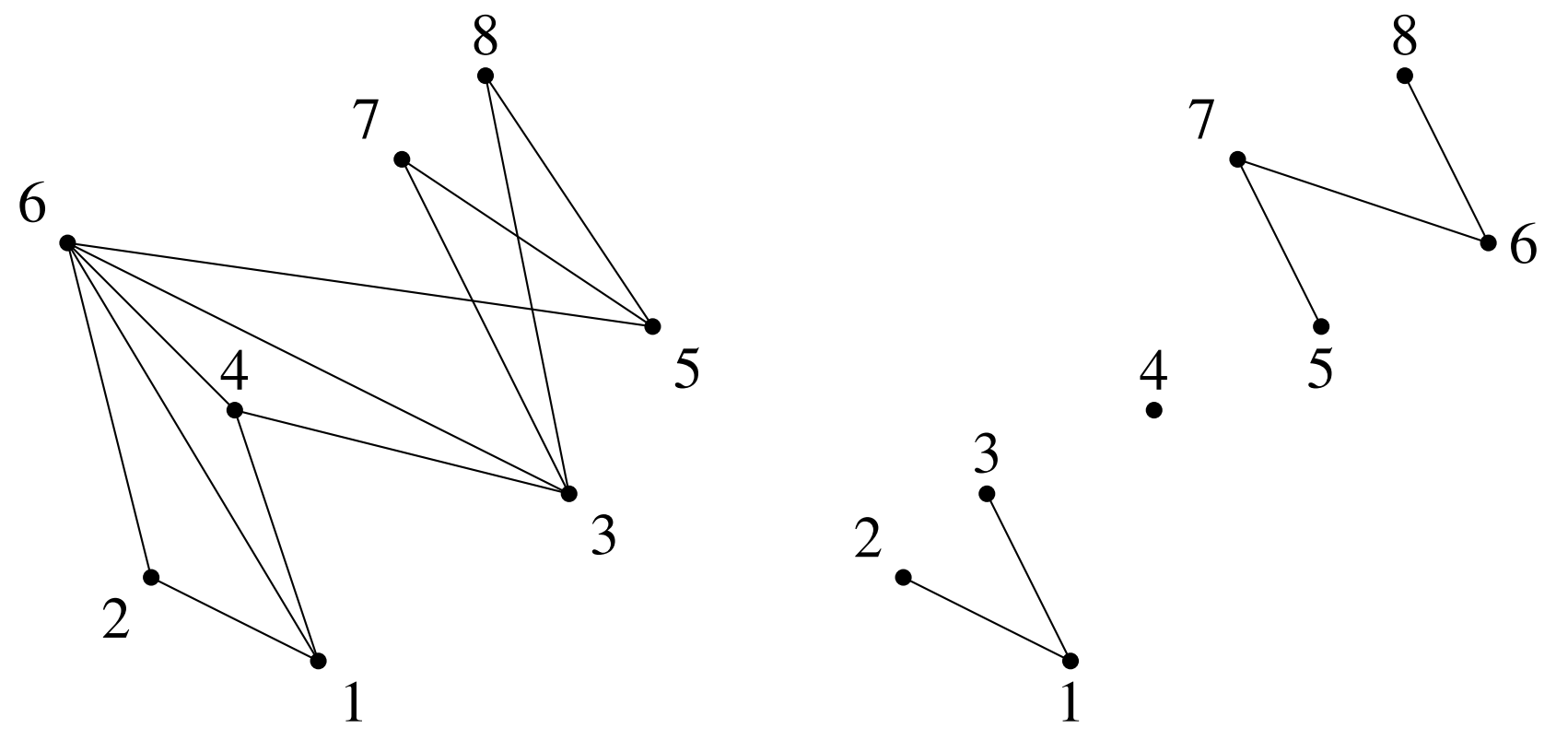

Let \(\pi \in S_n\). The inversion graph \(G_{\pi}\) has vertex set \(V(G_{\pi}) = [n]\) and edges
\begin{equation*}
  uv \in E(G_{\pi}) \iff (u - v)(\pi^{-1}(u) - \pi^{-1}(v)) < 0.
\end{equation*}
In other words, we think of the vertices as labeled by the \emph{entries} of \(\pi\), rather than the indices. In particular, we have that \(|E(G_{\pi})| = \inv(\pi)\), where the latter denotes the number of inversions of \(\pi\). The graph-theoretical notation we will use is standard, see e.g.\ Diestel \cite{diestel_graph_2025}. A convenient way to visualize \(G_{\pi}\) is to draw the plot of \(\pi\) in the Cartesian plane and connect any two points such that one is southeast of the other by an edge (see Figure~\ref{fig:inv graph example}).

Our permutation \(\pi\) is said to contain the pattern \(\sigma \in S_m\) (\(m \leq n\)) if \(\pi\) has a subsequence with the same relative order as \(\sigma\), both viewed in one-line notation. Otherwise, \(\pi\) avoids \(\sigma\). If \(B\) is a collection of patterns, we write \(\Av(B)\) for the set of all permutations avoiding every pattern in \(B\), and \(\Av_n(B) = \Av(B) \cap S_n\). We also write \(\Sub(B)\) for the set of all patterns contained in some \(\pi \in B\). Both are examples of \emph{permutation classes}, i.e.\ down-sets in the poset of all permutations ordered by pattern containment. There are extensive surveys on this topic by Bóna \cite{bona_combinatorics_2022}, Kitaev \cite{kitaev_patterns_2011} and Vatter \cite{vatter_permutation_2015}.

Let \(I \subseteq [n]\), and write \(\pi_I\) for the subsequence of \(\pi\) indexed by \(I\). Let \(\sigma\) be the pattern having the same relative order as \(\pi_I\). Then it is clear that the induced subgraph \(G_{\pi}[\pi(I)]\) is isomorphic to \(G_{\sigma}\) under the order-preserving map sending the set of entries \(\pi(I)\) to \([|I|]\). In this sense, we may say that patterns in \(\pi\) correspond with induced subgraphs of \(G_{\pi}\). For example, the subsequence \(643\) in Figure~\ref{fig:inv graph example} (left) has induced subgraph \(K_3\), and this is the same as the inversion graph of the corresponding pattern \(321\).

The direct sum of two permutations \(\pi \in S_n\) and \(\tau \in S_{\ell}\) is the permutation \(\pi \oplus \tau \in S_{n + \ell}\) defined by
\begin{equation*}
  (\pi \oplus \tau)(i) = \begin{cases}
    \pi(i) & \text{when } i \leq n, \\
    \tau(i - n) + n & \text{when } i > n.
  \end{cases}
\end{equation*}
A permutation is called indecomposable if it cannot be written as a direct sum of two nonempty permutations. Any permutation \(\pi\) decomposes uniquely into a direct sum of indecomposable permutations, called the components of \(\pi\). A related notion is the skew sum
\begin{equation*}
  (\pi \ominus \tau)(i) = \begin{cases}
    \pi(i) + \ell & \text{when } i \leq n, \\
    \tau(i - n) & \text{when } i > n.
  \end{cases}
\end{equation*}

We have that \(G_{\pi \oplus \tau} \cong G_{\pi} \sqcup G_{\tau}\) under the bijection \([n + \ell] \to [n] \sqcup [\ell]\) fixing the first \(n\) numbers and sending the last \(\ell\) to a disjoint copy of \([\ell]\). Thus, it is clear that the components of \(\pi\) correspond with the connected components of \(G_{\pi}\). An example of this is in Figure~\ref{fig:inv graph example} (right): 
\begin{equation*}
  \pi = 23147586 = 231 \oplus 1 \oplus 3142,
\end{equation*}
and the components of \(G_{\pi}\) are \(G_{231}\), \(G_1\) and \(G_{3142}\).

We will often run into permutations \(\pi\) that are obtained by inserting entries into another permutation \(\tau\) (so that \(\tau\) is contained in \(\pi\)). We want to visualize this process using the plot of \(\pi\), but it is convenient to relax the restriction that points are drawn exactly on the coordinates \(\{(i, \pi(i)) : i \in [n]\}\). Instead, we say that a \emph{diagram} of \(G_{\pi}\) is a configuration of \(n\) points in the plane so that the sequence of \(y\)-coordinates read from left to right has the same relative order as \(\pi\). As before, two points are connected by an edge if one is southeast of the other. We can then obtain a diagram of \(G_{\pi}\) by inserting points into a diagram of \(G_{\tau}\) without changing any of the previous points.

\section{Induced ears of 2-connected graphs} \label{sec:2-connected}

A classical tool for studying \(2\)-connected graphs is Whitney's open ear decomposition theorem \cite{whitney_non-separable_1932}. The statement, in its most common form, is that any \(2\)-connected graph \(G\) admits a sequence of subgraphs \(H_0 \subseteq H_1 \subseteq \ldots \subseteq H_r = G\) such that \(H_0\) is a cycle, and \(H_{i+1} = H_i \cup P\) for some \emph{\(H_i\)-path} \(P\). An \(H\)-path is a path in \(G\) with distinct endpoints in \(H\) and internal vertices in \(G \setminus H\). If \(P\) consists of a single edge, the edge is not in \(H\).

We would like to use this theorem for inductive arguments on inversion graphs. In our setting, however, the problem is that if the \(H_i\)-path is an edge, then \(H_i\) is not an induced subgraph of \(H_{i+1}\). In other words, if \(G = G_{\pi}\) is an inversion graph, we cannot guarantee that every \(H_i\) is an inversion graph. We will, therefore, need a modified version of the statement. It was difficult to find our version stated explicitly in the literature. It is based on the following lemma, which follows immediately from a stronger result of Kriesell \cite[Lemma~3]{kriesell_contractible_2000}.

\begin{lemma} \label{lem:ear}
  Let \(H\) be a proper induced \(2\)-connected subgraph of a \(2\)-connected graph \(G\). There exists a shortest \(H\)-path \(P\). It satisfies \(|E(P)| \geq 2\), and \(H \cup P\) is \(2\)-connected. If \(|E(P)| \geq 3\), then \(H \cup P\) is induced. 
\end{lemma}

\begin{proof}
  Let \(K\) be a component of \(G \setminus H\). Since \(G\) is \(2\)-connected, \(K\) has at least two distinct neighbors in \(H\). Hence, there exists an \(H\)-path with all internal vertices in \(K\). Let \(P\) be a shortest \(H\)-path. Observe that there is no \(H\)-path consisting of a single edge, since \(H\) is induced. Adding an \(H\)-path preserves \(2\)-connectivity, so \(H \cup P\) is \(2\)-connected. 
  
  Suppose \(P\) has at least two internal vertices. Denote its endpoints by \(x\) and \(y\). If \(uv \notin E(P)\) is an edge between two internal vertices of \(P\), then the segment \(uPv\) of \(P\) can be replaced by the shorter segment \(uv\). Similarly, if \(uv \notin E(P)\) is an edge from an internal vertex of \(P\) to a vertex of \(H\), then either \(xPuv\) or \(vuPy\) is a shorter \(H\)-path. Hence, \(H \cup P\) is induced.
\end{proof}

\begin{definition}
  Let \(G\) be a graph. An \emph{induced ear decomposition} of \(G\) is a sequence \(H_0 \subseteq H_1 \subseteq \ldots \subseteq H_r = G\) of induced subgraphs, where \(H_0\) is an induced cycle, and for each \(i\), either \(H_{i+1} = G[H_i \cup v]\) for a vertex \(v\) with at least two neighbors in \(H_i\), or \(H_{i+1} = H_i \cup P\) for some \(H_i\)-path \(P\).
\end{definition}

The existence of an induced ear decomposition for a \(2\)-connected graph is guaranteed by Lemma~\ref{lem:ear}.

\begin{lemma}
  A graph is \(2\)-connected if and only if it has an induced ear decomposition.
\end{lemma}

This is sufficient for the \(k = 2\) case of Theorem~\ref{thm:k connected}, but Theorem~\ref{thm:k connected av} requires something slightly more particular.

\begin{definition}
  Let \(G\) be a graph. A \emph{greedy induced ear decomposition} of \(G\) is a sequence of subgraphs \(H_0 \subseteq H_1 \subseteq \ldots \subseteq H_r = G\) in which \(H_0\) is an induced cycle, and the remaining graphs are constructed inductively as follows. Assume \(H_0, \ldots, H_{i-1}\) are constructed and \(H_{i-1} \neq G\).
  \begin{enumerate}[(i)]
    \item If there is a vertex in \(V(G \setminus H_{i-1})\) with at least two neighbors in \(H_{i-1}\), then choose such a vertex \(v\) and let \(H_i = G[H_{i-1} \cup v]\).
    \item Otherwise, choose a shortest \(H_{i-1}\)-path \(P\) in \(G\) and let \(H_i = H_{i-1} \cup P\).
  \end{enumerate}
  Terminate when neither choice can be made.
\end{definition}

The following properties are elementary.

\begin{lemma}
  A graph is \(2\)-connected if and only if it has a greedy induced ear decomposition. Furthermore, every greedy induced ear decomposition is an induced ear decomposition.
\end{lemma}

If \(H_0, \ldots, H_r\) is a (greedy) induced ear decomposition of an inversion graph \(G_{\pi}\), then we can identify \(H_i\) with the inversion graph \(G_{\tau}\) of the pattern \(\tau\) in \(\pi\) induced by the vertex set of \(H_i\). This is why it is crucial that each \(H_i\) be an \emph{induced} subgraph of \(G_{\pi}\).

Recall that inversion graphs are known to contain no induced cycles of length five or longer; they are \emph{\(4\)-chordal} -- see e.g.\ Tenner and Vatter \cite[Proposition~2.3]{tenner_cyclomatic_2026}. In general, a graph \(G\) is called \emph{\(\ell\)-chordal} if it has no induced cycles of length at least \(\ell + 1\). This limits the lengths of the shortest \(H\)-paths. We will assume throughout that \(\ell \geq 3\).

\begin{lemma}
  Let \(G\) be a \(2\)-connected \(\ell\)-chordal graph. If \(H\) is a connected induced subgraph of \(G\) with \(|V(H)| \geq 2\) and \(P\) is an \(H\)-path such that \(H \cup P\) is induced, then \(2 \leq |E(P)| \leq \ell - 1\).
\end{lemma}

\begin{proof}
  Let \(u\) and \(v\) be the endpoints of \(P\), and let \(Q\) be a shortest path from \(u\) to \(v\) in \(H\). We claim that \(C = uPvQu\) is an induced cycle in \(G\). Indeed, \(P\) has at least one internal vertex since \(H\) is induced, so \(P\) and \(Q\) are independent and \(C\) is a cycle. Since \(Q\) is shortest, \(C\) is induced in \(H \cup P\), and thus in \(G\). It follows that \(C\) is of length at most \(\ell\), so \(P\) can have at most \(\ell - 1\) edges.
\end{proof}

\begin{corollary} \label{cor:ear length}
  Let \(H_0, \ldots, H_r\) be an induced ear decomposition of a \(2\)-connected \(\ell\)-chordal graph \(G\). For every \(i\), we have that
  \begin{equation*}
    |V(H_{i+1})| - |V(H_i)| \leq \ell - 2.
  \end{equation*}
\end{corollary}

In fact, this gives us a strengthening of the \(k = 2\) case of Theorem~\ref{thm:k connected}: it works for all \(4\)-chordal graphs, and generalizes to \(\ell\)-chordal graphs.

\begin{proposition} \label{prop:l chordal bound}
  If \(G\) is a \(2\)-connected \(\ell\)-chordal graph on \(n\) vertices, then
  \begin{equation*}
    |E(G)| \geq \frac{\ell - 1}{\ell - 2} n - \frac{\ell}{\ell - 2}.
  \end{equation*}
\end{proposition}

\begin{remark}
  \(2\)-connected inversion graphs satisfy the bound with \(\ell = 4\):
  \begin{equation*}
    |E(G)| \geq \frac{3}{2}n - 2.
  \end{equation*}
\end{remark}

\begin{proof}[Proof of Proposition~\ref{prop:l chordal bound}]
We will prove the claim by induction on \(n\). First, suppose that \(G\) is a cycle. Since \(G\) is \(\ell\)-chordal, it must be the case that \(n \leq \ell\). We have
  \begin{equation*}
    \frac{\ell - 1}{\ell - 2} n - \frac{\ell}{\ell - 2} = n + \frac{n - \ell}{\ell - 2} \leq n = |E(G)|,
  \end{equation*}
  so the claim holds.

  Now, suppose that \(G\) is not a cycle, and let \(H\) be the penultimate graph in an induced ear decomposition of \(G\). Denote \(m = |V(H)|\). Corollary~\ref{cor:ear length} shows that \(n - m \leq \ell - 2\), and obviously \(|E(G)| - |E(H)| \geq n - m + 1\). Since \(H\) is an induced subgraph of \(G\), it is also \(\ell\)-chordal. By the induction hypothesis, we have
  \begin{equation*}
    |E(H)| \geq \frac{\ell - 1}{\ell - 2} m - \frac{\ell}{\ell - 2},
  \end{equation*}
  and it follows that
  \begin{align*}
    |E(G)| &\geq |E(H)| + n - m + 1\\
    &\geq \frac{\ell - 1}{\ell - 2} m - \frac{\ell}{\ell - 2} + n - m + 1 \\
    &= \frac{\ell - 1}{\ell - 2} n - \frac{\ell}{\ell - 2} + \frac{\ell - 2 - (n-m)}{\ell - 2}.
  \end{align*}
  The numerator of the last term was established to be nonnegative. Hence,
  \begin{equation*}
    |E(G)| \geq \frac{\ell - 1}{\ell - 2} n - \frac{\ell}{\ell - 2}. \qedhere
  \end{equation*}
\end{proof}

\section{Patterns in 2-connected inversion graphs} \label{sec:2-connected av}

This section builds up to a proof of the following \(k = 2\) version of Theorem~\ref{thm:k connected av}.

\begin{theorem} \label{thm:2 connected av}
  Let \(\sigma\) be a pattern. The following are equivalent:
  \begin{enumerate}[(i)]
    \item\label{thm:2 connected av:forcing} There exists a constant \(C\) depending only on \(\sigma\) such that for every \(\pi \in \Av(\sigma)\) whose inversion graph is \(2\)-connected,
    \begin{equation*}
      \inv(\pi) \geq 2|\pi| - C.
    \end{equation*}
    \item\label{thm:2 connected av:patterns} We have \(\sigma \in \Av(321, 2341, 4123)\).
    \item\label{thm:2 connected av:graphs} The components of \(G_{\sigma}\) are paths and \(4\)-cycles.
  \end{enumerate}
  If \(\sigma\) satisfies these properties and \(G_{\sigma}\) has \(p_1\) isolated vertices, \(p_2\) components consisting of a single edge and \(q\) \(4\)-cycles, then we may choose
  \begin{equation*}
    C = |\sigma| + 2p_1 + p_2 - q.
  \end{equation*}
\end{theorem}

In the previous section, we could work purely in an abstract graph-theoretical setting. Those foundations will be useful for us here, but a proof of Theorem~\ref{thm:2 connected av} requires further structural analysis that is specific to inversion graphs.

\begin{definition}
  Let \(G_{\pi}\) be a \(2\)-connected inversion graph. An \emph{anchored ear decomposition} of \(G_{\pi}\) is a greedy induced ear decomposition \(H_0, \ldots, H_r\) such that the induced cycle \(H_0\) contains the vertices \(1\) and \(\pi(1)\).
\end{definition}

\begin{lemma}
  Every \(2\)-connected inversion graph has an anchored ear decomposition.
\end{lemma}

\begin{proof}
  It suffices to show that there exists an induced cycle containing the desired vertices. Let \(G_{\pi}\) be a \(2\)-connected inversion graph, and let \(i = \pi^{-1}(1)\). Observe first that we must have \(\pi(1) > 2\) (and, similarly, \(i > 2\)), as the degree of \(\pi(1)\) (resp.\ \(1\)) would otherwise be at most one.

  If there exists some \(j\) such that \(1 < j < i\) and \(1 < \pi(j) < \pi(1)\), then clearly \(\pi(1), \pi(j), 1\) form a \(3\)-cycle. Otherwise, we must have \(\pi(2) > \pi(1)\) and \(\pi^{-1}(2) > i\), so that the vertices \(\pi(1), \pi(2), 1, 2\) form a \(3412\)-pattern, i.e.\ an induced \(4\)-cycle. 
\end{proof}

Recall that by Corollary~\ref{cor:ear length}, if \(H_0, \ldots, H_r\) is an induced ear decomposition of a \(4\)-chordal graph \(G\), then \(H_i\) and \(H_{i+1}\) can differ by at most two vertices. If they differ by exactly two, then \(H_{i+1}\) is of the form \(H_i \cup P\) for some \(H_i\)-path \(P\) with two internal vertices. Our key observation is that if \(G\) is an inversion graph, then the way \(P\) attaches to \(H_i\) is highly constrained.

If \(H_0, \ldots, H_r\) is an induced ear decomposition of a \(2\)-connected inversion graph \(G_{\pi}\), we can build a diagram of \(G_{\pi}\) by starting from a diagram of \(H_0\), and inductively inserting points so that the \(i\)-th diagram represents \(H_i\). We do not need to move the pre-existing points of the diagram when the new points are inserted. Anchored ear decompositions have the distinctive property that the leftmost and bottommost points of the diagrams stay fixed throughout.

\begin{figure}[t]
  \centering
  \begin{tikzpicture}[scale=0.46,
      every label/.style={font=\scriptsize, inner sep=1pt},
      label distance=0.5pt,
      v1/.style={label=above:\(3\)},
      v2/.style={label=left:\(7\)},
      v3/.style={label=above right:\(4\)},
      v4/.style={label=right:\(1\)},
      v5/.style={label=above right:\(5\)},
      v6/.style={label=left:\(8\)},
      v7/.style={label=below:\(2\)},
      v8/.style={label=below:\(6\)}]
    \foreach \x\y in {1/3, 3/4, 4/1, 7/2}
      \coordinate (\x) at (\x,\y);
    \draw[rpath, thick]
      (1) -- (4) -- (3) -- (7) -- cycle;
    \foreach \x in {1, 3, 4, 7}
      \node[dot, v\x] at (\x) {};

    \node at (4,-0.3) {\(H_0\)};

    \begin{scope}[shift={(7,0)}]
      \foreach \x\y in {1/3, 2/7, 3/4, 4/1, 7/2}
      \coordinate (\x) at (\x,\y);
      \draw[rpath, very thin]
        (1) -- (4) -- (3) -- (7) -- cycle;
      \foreach \x in {1, 3, 4, 7}
        \node[dot, v\x] at (\x) {};

      \draw[rpath, thick]
        (2) -- (3)
        (2) to[bend right=20] (4)
        (2) -- (7);
      \node[dot, v2] at (2) {};
      \node at (4,-0.3) {\(H_1\)};
    \end{scope}

    \begin{scope}[shift={(14,0)}]
      \foreach \x\y in {1/3, 2/7, 3/4, 4/1, 6/8, 7/2, 8/6}
      \coordinate (\x) at (\x,\y);
      \draw[rpath, very thin]
        (1) -- (4) -- (3) -- (7) -- cycle
        (2) -- (3)
        (2) to[bend right=20] (4)
        (2) -- (7);
      \foreach \x in {1, 2, 3, 4, 7}
        \node[dot, v\x] at (\x) {};

      \draw[rpath, thick]
        (7) -- (6) -- (8) -- (2);
      \foreach \x in {6,8}
        \node[dot, v\x] at (\x) {};
      \node at (4.5,-0.3) {\(H_2\)};
    \end{scope}

    \begin{scope}[shift={(22,0)}]
      \foreach \x\y in {1/3, 2/7, 3/4, 4/1, 5/5, 6/8, 7/2, 8/6}
      \coordinate (\x) at (\x,\y);
      \draw[rpath, very thin]
        (1) -- (4) -- (3) -- (7) -- cycle
        (2) -- (3)
        (2) to[bend right=20] (4)
        (2) -- (7)
        (7) -- (6) -- (8) -- (2);
      \foreach \x in {1, 2, 3, 4, 6, 7, 8}
        \node[dot, v\x] at (\x) {};

      \draw[rpath, thick]
        (2) -- (5) -- (7);
      \node[dot, v5] at (5) {};
      \node at (4.5,-0.3) {\(H_3\)};
    \end{scope}
  \end{tikzpicture}
  \caption{Diagrams of an induced ear decomposition \(H_0, H_1, H_2, H_3\) of the inversion graph of \(37415826\).} 
  \label{fig:ears}
\end{figure}
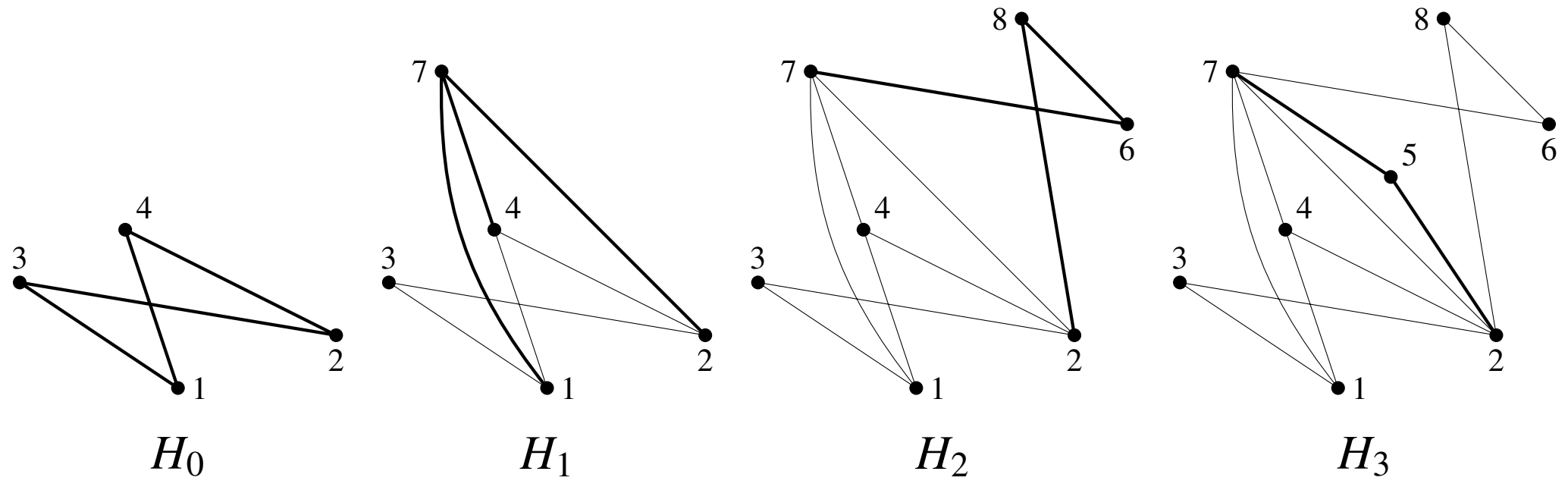

\begin{example}
  Let \(\pi = 37415826\), and consider the induced ear decomposition \(H_0, H_1, H_2, H_3\), where \(H_0\) has vertex set \(\{1, 2, 3, 4\}\), \(H_1\) has vertex set \(\{1, 2, 3, 4, 7\}\), \(H_2\) has vertex set \(\{1, 2, 3, 4, 6, 7, 8\}\) and \(H_3 = G_{\pi}\); see Figure~\ref{fig:ears}. As per the figure, a diagram for each graph can be obtained by adding points and edges to the previous diagram. This ear decomposition is not greedy, since vertex \(5\) could have been included in \(H_2\) instead of \(\{6, 8\}\). If we make this change, the resulting ear decomposition is anchored, since the first graph contains \(1\) and \(\pi(1) = 3\).

\end{example}

\begin{figure}[t]
  \centering
  \begin{tikzpicture}[scale=0.7]
    \draw[fill=black!20] (0,0) rectangle (5,5);
    \foreach \x\y\l in {2.6/6/y, 6/2.3/x, 5.5/7/u, 7/5.5/v}
      \coordinate (\l) at (\x,\y);

    \draw[rpath]
      (8,5) -- (5,5) -- (5,8);
    \draw[rpath, very thin]
      (x) -- (y)
      (y) to[bend right=8] (3.1,3.4)
      (y) to[bend right=5] (3.6,4.0)
      (x) to[bend left=10] (3.4,2.85)
      (x) to[bend left=5] (4.0,2.9);
    \draw[rpath, thick]
      (x) -- (u) -- (v) -- (y);

    \node at (2.5,2.5) {\(H_i\)};
    \foreach \l in {x, y, u, v}
      \node[dot] at (\l) {};
    \node[above] at (u) {\(u\)};
    \node[right] at (v) {\(v\)};
    \node at (4,-1) {(i)};

    \begin{scope}[shift={(12,3)}]
      \draw[fill=black!20] (0,0) rectangle (5,5);
      \foreach \x\y\l in {-1/2.7/y, 2.4/-1/x, -2/-0.5/u, -0.5/-2/v}
        \coordinate (\l) at (\x,\y);

      \draw[rpath]
        (-3,0) -- (0,0) -- (0,-3);
      \draw[rpath, very thin]
        (x) -- (y)
        (x) to[bend right=8] (1.9,1.6)
        (x) to[bend right=5] (1.4,1.0)
        (y) to[bend left=10] (1.6,2.15)
        (y) to[bend left=5] (1.0,2.1);
      \draw[rpath, thick]
        (x) -- (u) -- (v) -- (y);

      \node at (2.5,2.5) {\(H_i\)};
      \foreach \l in {x, y, u, v}
        \node[dot] at (\l) {};
      \node[left] at (u) {\(u\)};
      \node[below] at (v) {\(v\)};
      \node at (1,-4) {(ii)};
    \end{scope}
  \end{tikzpicture}
  \caption{The two options in Lemma~\ref{lem:inv graph 2 ears}. In (i), the shaded square contains every point of \(H_i\) except its topmost and rightmost points; in (ii), every point except its bottommost and leftmost points.}
  \label{fig:2 ears}
\end{figure}
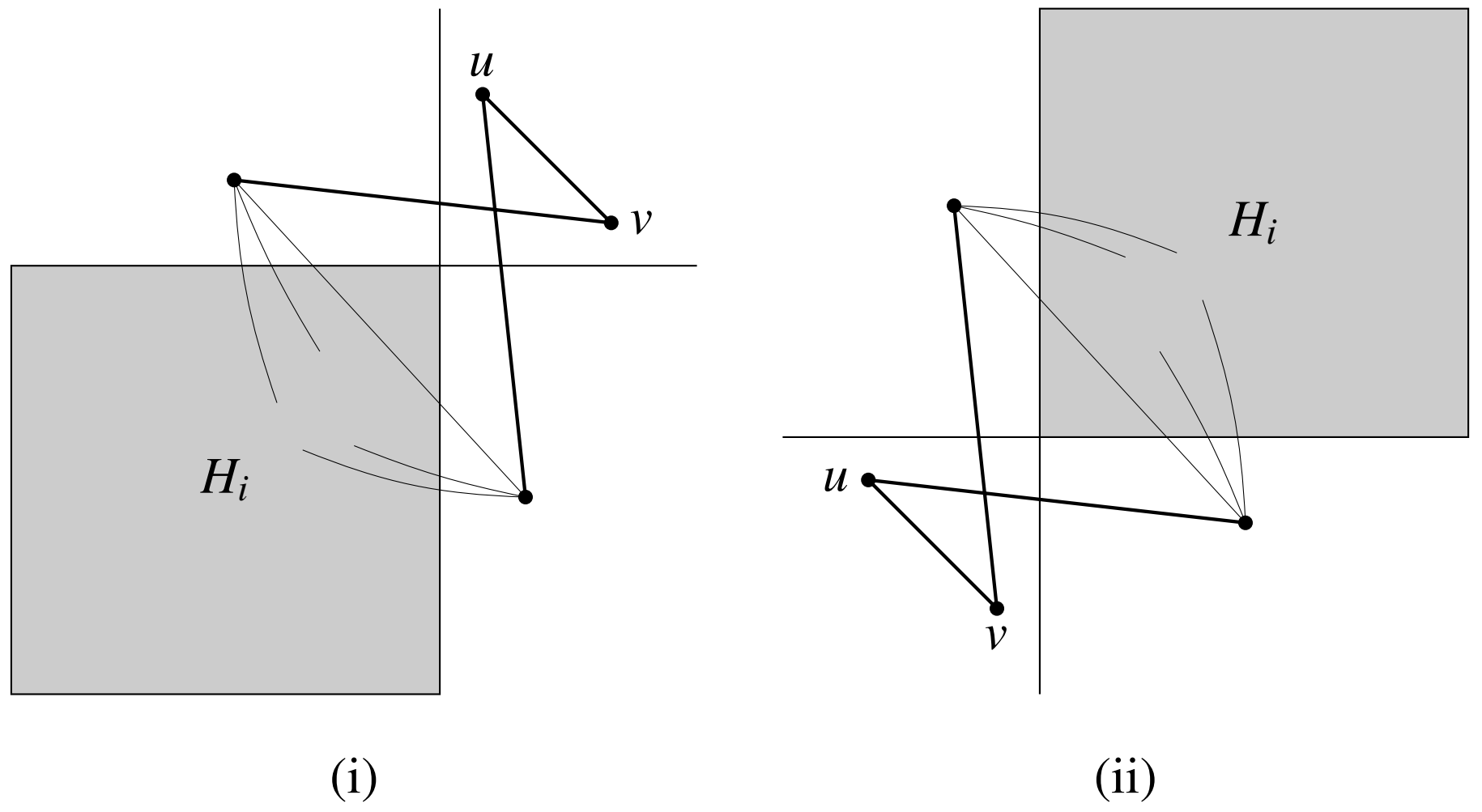

\begin{lemma} \label{lem:inv graph 2 ears}
  Let \(H_0, \ldots, H_r\) be an induced ear decomposition of a \(2\)-connected inversion graph, and suppose that \(H_i\) and \(H_{i+1}\) differ by two vertices \(u\) and \(v\) with \(u > v\). Let \(D_i\) be a diagram of \(H_i\), and \(D_{i+1}\) a diagram of \(H_{i+1}\) that is obtained by inserting \(u\) and \(v\) into \(D_i\). Then \(D_{i+1}\) satisfies one of the following:
  \begin{enumerate}[(i)]
    \item \label{lem:inv graph 2 ears:top right} 
    \(u\) is topmost and third from the right, and \(v\) is rightmost and third from the top; or
    \item \label{lem:inv graph 2 ears:bot left} \(u\) is leftmost and third from the bottom, and \(v\) is bottommost and third from the left.
  \end{enumerate}
  If the decomposition is anchored, then only option \ref{lem:inv graph 2 ears:top right} is possible. Figure~\ref{fig:2 ears} visualizes both options.

\end{lemma}

\begin{remark}
  This lemma shows that for even \(n \geq 4\), there exists a unique permutation \(\pi \in S_n\) with the minimum number of inversions whose inversion graph is \(2\)-connected. This permutation is determined by its induced ear decomposition \(H_0, \ldots, H_r\), in which \(H_0\) is a \(4\)-cycle, and each step adds two vertices. By the lemma, \(G_{\pi}\) must be the ladder graph; see Figure~\ref{fig:pi big} in the appendix. The first few such permutations are
  \begin{equation*}
    3412,\quad 351624,\quad 35172846.
  \end{equation*}
\end{remark}

\begin{proof}[Proof of Lemma~\ref{lem:inv graph 2 ears}]
  Denote the unique neighbors of \(u\) and \(v\) in \(H_i\) by \(x\) and \(y\), respectively. In \(H_{i+1} \setminus x\), every remaining vertex of \(H_i\) must lie southwest or northeast of \(u\). Since \(H_i\) is \(2\)-connected, they must all lie on the \emph{same} side; assume they all lie southwest of \(u\). Symmetrically, every point of \(H_i \setminus y\) must lie southwest or northeast of \(v\). If they are northeast of \(v\), then \(u\) also lies northeast of \(v\). This is impossible, since \(uv\) is an edge.

  Hence, every point of \(H_{i} \setminus \{x, y\}\) lies southwest of both \(u\) and \(v\). Since \(x\) and \(y\) are distinct, it is easy to see that condition \ref{lem:inv graph 2 ears:top right} is satisfied. If, instead, every point of \(H_{i} \setminus \{x, y\}\) lies northeast of both \(u\) and \(v\), then condition \ref{lem:inv graph 2 ears:bot left} holds.
\end{proof}

\begin{figure}[t]
  \centering
  \begin{tikzpicture}[scale=0.7]
    \draw[fill=black!20] (0,0) rectangle (5,5);
    \fill[black!20] (5,5) rectangle (10,10);
    \foreach \x\y\l in {2.6/6/y, 6/2.3/x, 5.5/7/u, 7/5.5/v, 2/5/a, 5/3/b} 
      \coordinate (\l) at (\x,\y);

    \draw[rpath]
      (0,5) -- (10.2,5)
      (5,0) -- (5,10.2);
    \draw[rpath, very thin]
      (x) -- (y)
      (y) to[bend right=8] (3.1,3.4)
      (y) to[bend right=5] (3.6,4.0)
      (x) to[bend left=10] (3.4,2.85)
      (x) to[bend left=5] (4.0,2.9);
    \draw[rpath, thick]
      (x) -- (u) -- (v) -- (y);

    \node at (2.5,2.5) {\(H_i\)};
    \node at (2.5,7.5) {\(Y\)};
    \node at (7.5,2.5) {\(X\)};
    \foreach \l in {x, y, u, v, a, b}
      \node[dot] at (\l) {};
    \node[above] at (u) {\(u\)};
    \node[right] at (v) {\(v\)};
    \node[right] at (x) {\(x\)};
    \node[above] at (y) {\(y\)};
  \end{tikzpicture}
  \caption{Regions \(X\) and \(Y\) in an anchored ear decomposition containing no other points.}
  \label{fig:empty X Y}
\end{figure}
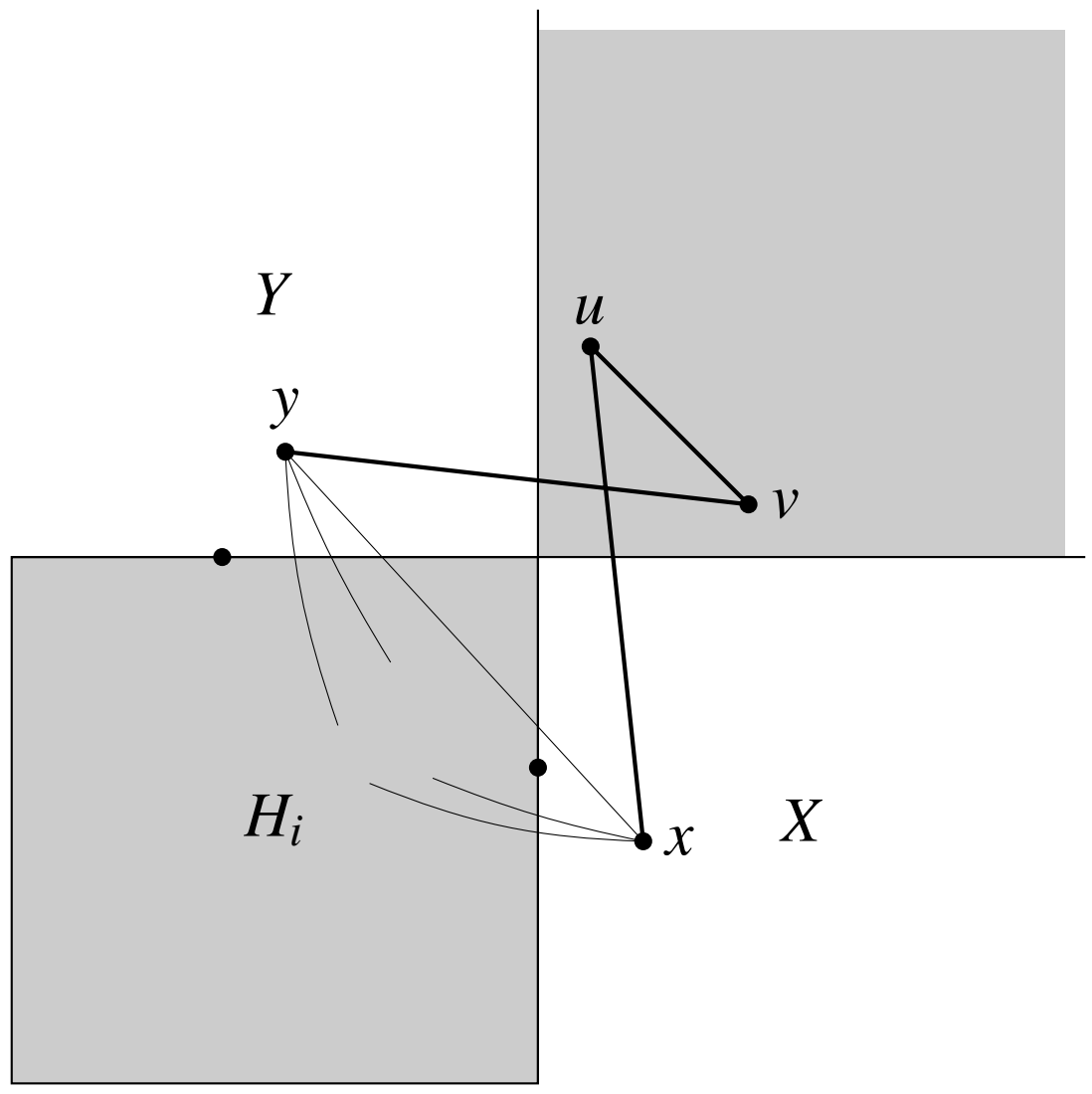

The convenience of anchored ear decompositions lies precisely in the fact that only the former of the two options of Lemma~\ref{lem:inv graph 2 ears}, i.e.\ Figure~\ref{fig:2 ears} (i), is possible. Indeed, in any anchored ear decomposition of a \(2\)-connected inversion graph, we know exactly how every \(2\)-vertex ear must be inserted. The next lemma will explain why we insist particularly upon \emph{greedy} ear decompositions.

Let \(H_0, \ldots, H_r\) be an anchored ear decomposition of a \(2\)-connected inversion graph \(G_{\pi}\), and suppose that \(H_i\) and \(H_{i+1}\) differ by two vertices \(u > v\). Denote by \(x\) and \(y\) their unique neighbors in \(H_i\). In a diagram \(D\) of \(G_{\pi}\), draw a vertical line through the point of \(H_i\) that is immediately left of \(x\), and a horizontal line through the point of \(H_i\) that is immediately below \(y\), as in Figure~\ref{fig:empty X Y}. Let \(Y\) and \(X\), respectively, denote the northwest and southeast regions of \(D\) cut out by these lines.

\begin{lemma} \label{lem:empty X Y}
  The interiors of the regions \(X\) and \(Y\) defined above (see Figure~\ref{fig:empty X Y}) contain no points of \(D\) other than \(x\) and \(y\), respectively.
\end{lemma}

\begin{proof}
  Suppose that \(X\) or \(Y\) contains some other point \(z\). This point must have been added in some step \(H_j \to H_{j+1}\) with \(j > i\) (recall that \(u\) and \(v\) were added in the step \(H_i \to H_{i+1}\)). However, \(z\) has at least two edges to \(H_i\): if it lies in \(X\), then it has an edge to both \(y\) and the point immediately below it, and the other case is symmetrical. This contradicts the greediness of our ear decomposition, as we should have added \(z\) instead of \(\{u, v\}\) in the step \(H_i \to H_{i+1}\).
\end{proof}

We are now nearly ready to prove a central result that almost immediately implies the \(k=2\) case of Theorem~\ref{thm:k connected av}. This is, really, the main mechanism connecting pattern avoidance to our graph-theoretical machinery. First, some notation.

\begin{definition}
  We call a pattern \(\sigma\) \emph{\(k\)-forcing} if there exists a constant \(C\) such that for every permutation \(\pi \in \Av(\sigma)\) whose inversion graph is \(k\)-connected, we have that
  \begin{equation} \label{eq:k forcing}
    \inv(\pi) \geq k|\pi| - C.
  \end{equation}
  If \(\sigma\) is \(k\)-forcing, let \(C_k(\sigma)\) denote the smallest constant \(C \geq 3\) satisfying \eqref{eq:k forcing}.
\end{definition}

Our goal is to show that the permutation class \(\mathcal C\) of \(2\)-forcing patterns (and later, the patterns that are \(k\)-forcing for every \(k \geq 2\)) is exactly \(\Av(321, 2341, 4123)\). This will follow from four separate observations: \(\mathcal C\) contains \(3412\); \(\mathcal C\) contains all permutations whose inversion graphs are paths; \(\mathcal C\) is closed under direct sums; and, lastly, the three patterns \(321\), \(2341\) and \(4123\) are not in \(\mathcal C\). The third of these is the main mechanism alluded to above.

For purposes of notational convenience, we define the \emph{deficit} of a graph \(G\) as
\begin{equation*}
  \deficit(G) = 2|V(G)| - |E(G)|. 
\end{equation*}
Observe that \(2\)-forcing is equivalent to the statement that \(\deficit(G_{\pi})\) is uniformly bounded from above over every \(\pi \in \Av(\sigma)\) with a \(2\)-connected inversion graph. The following lemma is immediate.

\begin{lemma} \label{lem:4 chordal deficit}
  Let \(H_0, \ldots, H_r\) be an induced ear decomposition of a \(4\)-chordal graph \(G\).
  \begin{enumerate}[(i)]
    \item If \(H_i \to H_{i+1}\) adds exactly one vertex, then \(\deficit(H_{i+1}) \leq \deficit(H_i)\).
    \item If \(H_i \to H_{i+1}\) adds two vertices, then \(\deficit(H_{i+1}) = \deficit(H_i) + 1\).
  \end{enumerate}
\end{lemma}

\begin{theorem} \label{thm:direct sum closure}
  If \(\sigma\) and \(\tau\) are \(2\)-forcing patterns, then \(\sigma \oplus \tau\) is \(2\)-forcing with
  \begin{equation*}
    C_2(\sigma \oplus \tau) \leq C_2(\sigma) + C_2(\tau).
  \end{equation*}
  In particular, the permutation class of \(2\)-forcing patterns is closed under direct sums.
\end{theorem}

\begin{proof}
  Let \(\pi\) be a \((\sigma \oplus \tau)\)-avoiding permutation whose inversion graph is \(2\)-connected. Let \(H_0, \ldots, H_r\) be an anchored ear decomposition of \(G_{\pi}\). If \(\pi\) avoids \(\sigma\) itself, then we have that \(\inv(\pi) \geq 2|\pi| - C_2(\sigma)\), and we are done. So, suppose that \(\pi\) contains \(\sigma\), and let \(i\) be minimal so that \(H_i\) contains \(\sigma\).

  Suppose first that there are at most two indices \(j\) such that \(j \geq i\) and the step \(H_j \to H_{j+1}\) adds two vertices. Since \(H_{i-1}\) (if it exists) avoids \(\sigma\), Lemma~\ref{lem:4 chordal deficit} implies that
  \begin{equation*}
    \deficit(G_{\pi}) \leq \deficit(H_i) + 2 \leq C_2(\sigma) + 3 \leq C_2(\sigma) + C_2(\tau).
  \end{equation*}
  Otherwise, let \(i \leq a < b < c\) be the three smallest indices for which the steps \(H_a \to H_{a+1}\), \(H_b \to H_{b+1}\) and \(H_c \to H_{c+1}\) add two vertices each. By Lemma~\ref{lem:empty X Y}, the diagram of \(G_{\pi}\) must be of the form described in Figure~\ref{fig:three ears}, with the interior of each white rectangle empty (except for the six attachment points of the ears); from southwest to northeast, the three thick paths are the ears added in steps \(a\), \(b\) and \(c\). Each remaining point of \(G_\pi\) must lie inside the shaded regions. 

  \begin{figure}[t]
    \centering
    \begin{tikzpicture}[scale=0.8]
      \draw[fill=black!20] (0,0) rectangle (4.5,4.5);
      \fill[black!20] (4.5,4.5) rectangle (8.5,8.5);
      \fill[black!20] (8.5,8.5) rectangle (12.5,12.5);
      \fill[black!20] (12.5,12.5) rectangle (16.5,16.5);

      \draw[rpath]
        (0,4.5) -- (16.7,4.5)
        (0,8.5) -- (16.7,8.5)
        (0,12.5) -- (16.7,12.5)
        (4.5,0) -- (4.5,16.7)
        (8.5,0) -- (8.5,16.7)
        (12.5,0) -- (12.5,16.7);

      \node at (2.25,2.25) {\(H_a\)};
      \node at (14.5,14.5) {\(R\)};

      \foreach \x\y\l in {
          2.5/5.35/ya, 5.45/2.1/xa, 4.9/6.1/ua, 6.1/4.9/va, 1.8/4.5/pa, 4.5/2.75/qa,
          6.75/9.15/yb, 9.2/6.55/xb, 8.9/10.1/ub, 10.1/8.9/vb, 6.4/8.5/pb, 8.5/7.05/qb,
          10.45/13.5/yc, 13.4/10.25/xc, 12.9/14.1/uc, 14.1/12.9/vc, 9.9/12.5/pc, 12.5/10.8/qc}
        \coordinate (\l) at (\x,\y);

      \draw[rpath, very thin]
        (xa) -- (ya)
        (ya) to[bend right=6] (2.75,3.05)
        (ya) to[bend right=10] (3.6,3.6)
        (xa) to[bend left=7] (2.95,2.35)
        (xa) to[bend left=4] (3.75,2.8)
        (xb) -- (yb)
        (yb) to[bend right=5] (7.1,7.55)
        (yb) to[bend right=10] (7.55,7.65)
        (xb) to[bend left=5] (6.6,7.15)
        (xb) to[bend left=9] (7.6,6.75)
        (xc) -- (yc)
        (yc) to[bend right=3] (10.55,10.8)
        (yc) to[bend right=8] (11.2,11.9)
        (xc) to[bend left=6] (11.0,10.5)
        (xc) to[bend left=10] (12.0,10.65);
      \draw[rpath, thick]
        (xa) -- (ua) -- (va) -- (ya)
        (xb) -- (ub) -- (vb) -- (yb)
        (xc) -- (uc) -- (vc) -- (yc);

      \foreach \l in {xa, ya, ua, va, pa, qa, xb, yb, ub, vb, pb, qb, xc, yc, uc, vc, pc, qc}
        \node[dot] at (\l) {};
      \node[right] at (xa) {\(x\)};
      \node[above] at (ya) {\(y\)};
      \node[right] at (xc) {\(w\)};
      \node[above] at (yc) {\(z\)};
    \end{tikzpicture}
    \caption{Three ears in the diagram of \(G_{\pi}\) for the proof of Theorem~\ref{thm:direct sum closure}.}
    \label{fig:three ears}
  \end{figure}
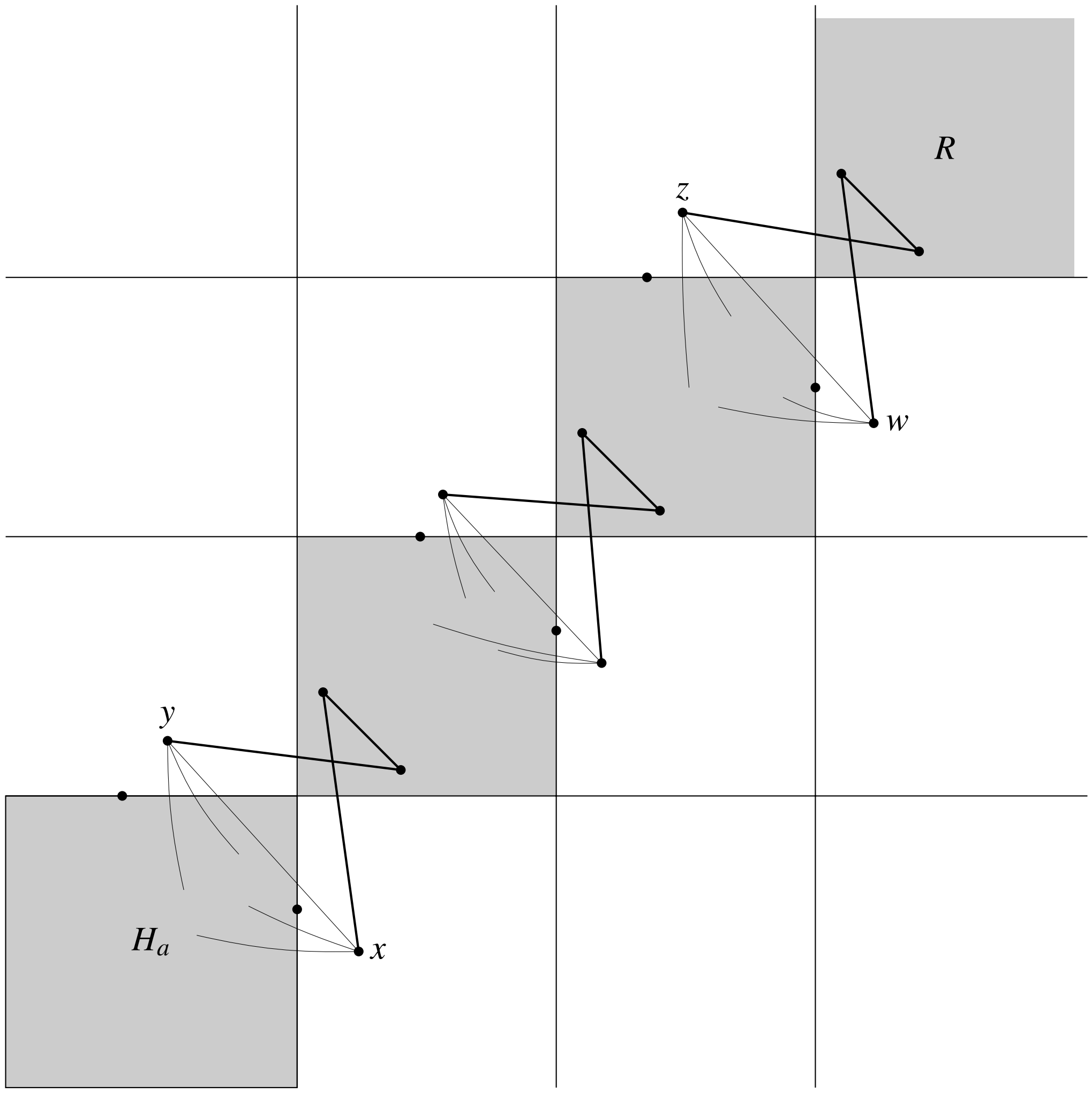

  Let \(x\) and \(y\) be the rightmost and topmost points of \(H_a\), and \(w\) and \(z\) the rightmost and topmost points of \(H_c\); these are the neighbors of the ears added in steps \(a\) and \(c\), respectively. We claim that \(z\) and \(w\) lie strictly northeast of both \(x\) and \(y\). Indeed, using Lemmas~\ref{lem:inv graph 2 ears} and \ref{lem:empty X Y}, \(x\) must lie weakly to the west of the middle vertical line of the figure, whereas \(z\) must lie strictly east of it. Similarly, \(y\) is weakly south of the middle horizontal line, and \(w\) is strictly north of it. Since \(x\) is southeast of \(y\) and \(w\) is southeast of \(z\), the claim follows.

  Consider now the subgraph \(R\) of \(G_{\pi}\) induced on the vertices in the northeastern shaded square of the figure together with \(z\) and \(w\). \(R\) is \(2\)-connected, since our induced ear decomposition of \(G_{\pi}\) restricts to an induced ear decomposition of \(R\) with initial cycle consisting of \(\{z, w\}\) and their neighbors added in step \(H_c \to H_{c+1}\). Furthermore, crucially, \(R\) must avoid \(\tau\) due to the discussion in the previous paragraph: if \(R\) contains \(\tau\), then this occurrence of \(\tau\) together with the occurrence of \(\sigma\) in \(H_{a}\) forms an occurrence of \(\sigma \oplus \tau\) in \(G_{\pi}\).

  Denote now by \(L\) the induced subgraph on vertex set \((V(G_{\pi}) \setminus V(R)) \cup \{z, w\}\), so that \(L \cup R = G_{\pi}\) with the two graphs intersecting exactly at the edge \(zw\). \(L\) is also a \(2\)-connected inversion graph, since the induced ear decomposition of \(G_{\pi}\) restricts to an induced ear decomposition of \(L\). By Lemma~\ref{lem:4 chordal deficit}, we have
  \begin{align*}
    \deficit(L) &\leq \deficit(H_a) + 2 \leq C_2(\sigma) + 3, \\*
    \deficit(R) &\leq C_2(\tau),
  \end{align*}
  and thus
  \begin{align*}
    \deficit(G_{\pi}) &= 2(|V(L)| + |V(R)| - 2) - (|E(L)| + |E(R)| - 1) \\
    &= \deficit(L) + \deficit(R) - 3 \\
    &\leq C_2(\sigma) + C_2(\tau).
  \end{align*}
  This concludes the proof.
\end{proof}

Let us now embark on the remaining, easier observations we need. The following two lemmas are noteworthy mainly because this is where \emph{specific} patterns come into play.

\begin{lemma} \label{lem:3412 is 2-forcing}
  The pattern \(3412\) is \(2\)-forcing with \(C_2(3412) = 3\).
\end{lemma}

\begin{proof}
  Since inversion graphs are \(4\)-chordal and \(3412\) is the only permutation whose inversion graph is a \(4\)-cycle, the inversion graphs of \(3412\)-avoiding permutations are chordal. The claim then follows from classical graph theory -- namely, Dirac's theorem \cite{dirac_rigid_1961} that every chordal graph has a simplicial vertex.

  Using our tools, we can also observe that if \(H_0, \ldots, H_r\) is an induced ear decomposition of a \(2\)-connected inversion graph and \(H_i \to H_{i+1}\) adds two new vertices, then \(H_{i+1}\) contains a \(3412\)-pattern by Lemma~\ref{lem:inv graph 2 ears} (see Figure~\ref{fig:2 ears}). Hence, if \(\pi \in \Av(3412)\), then
  \begin{equation*}
    \deficit(G_{\pi}) \leq \deficit(H_0) = \deficit(C_3) = 3
  \end{equation*}
  by Lemma~\ref{lem:4 chordal deficit}.
\end{proof}

\begin{lemma} \label{lem:paths are 2-forcing}
  Let \(\sigma\) be a pattern whose inversion graph is a path on \(m\) vertices. Then \(\sigma\) is \(2\)-forcing with \(C_2(\sigma) \leq \max\{3, m\}\).
\end{lemma}

\begin{proof}
  Let \(H_0, \ldots, H_r\) be an anchored ear decomposition of a \(2\)-connected inversion graph \(G_{\pi}\). For this proof, it is more convenient to assume that \(H_0\) consists of the single edge connecting \(\pi(1)\) and \(1\), and that \(H_0 \to H_1\) is a step adding one or two vertices, in either case yielding an induced cycle containing \(\{1, \pi(1)\}\). We have \(\deficit(H_0) = 3\), and Lemmas~\ref{lem:empty X Y} and \ref{lem:4 chordal deficit} hold verbatim. 
  
  Let \(i_1 < \ldots < i_t\) be the set of all the indices such that the step \(H_{i_j} \to H_{i_j + 1}\) adds two new vertices; note that if there are no such steps, then 
  \begin{equation*}
    \deficit(G_{\pi}) \leq \deficit(H_0) = 3,
  \end{equation*}
  and we are done. So, we may assume that \(t \geq 1\). For each \(j\), let \(B_j\) be the set consisting of the topmost and rightmost vertices of \(H_{i_j}\). 
  
  If \(m \leq 2\) the claim holds vacuously, so assume \(m \geq 3\). It is easy to see that there are exactly two patterns \(\sigma\) of size \(m \geq 3\) whose inversion graphs are paths. They are distinguished by whether the path starts at vertex \(1\) or \(\sigma(1)\), and they are inverses of each other; see Albert, Brignall, Ruškuc and Vatter \cite[Section~9]{albert_rationality_2019} for details. We may assume without loss of generality that our pattern is of the former kind.

  Let \(P\) be a shortest path connecting \(1\) and the topmost vertex \(n\) in \(G_{\pi}\). In particular, \(P\) is induced. Furthermore, by Lemma~\ref{lem:empty X Y}, it must intersect every set \(B_j\). These sets are disjoint and do not contain \(n\), so \(P\) has at least \(t + 1\) vertices.
  
  We can do slightly better: if \(P\) has \emph{exactly} \(t+1\) vertices, then it must consist of \(n\) and exactly one vertex from each \(B_j\). The neighbor of \(n\) in \(P\) must be the southeastern vertex of \(B_t\). The rightmost vertex \(\pi(n)\) is not in \(\bigcup_j B_j\), and its only possible neighbor in \(\bigcup_j B_j\) is the northwestern vertex of \(B_t\). Thus, the path \(P\pi(n)\) is induced and contains \(t + 2\) vertices.
  
  In either case, \(G_{\pi}\) must contain an induced path \(Q\) beginning at \(1\) with at least \(t+2\) vertices. If \(t+2 \geq m\), then \(Q\) contains the pattern \(\sigma\). If \(\pi \in \Av(\sigma)\), we must therefore have \(t \leq m - 3\), i.e.\
  \begin{equation*}
      \deficit(G_{\pi}) \leq \deficit(H_0) + m - 3 = m. \qedhere
  \end{equation*}
\end{proof}

\begin{remark}
  Menger's theorem can be used for a different (and potentially shorter) proof, but we choose to flaunt our toolbox instead.
\end{remark}

With these tools, we are ready to prove the deepest part of Theorem~\ref{thm:2 connected av}, which is that its condition \ref{thm:2 connected av:patterns} (or, equivalently, \ref{thm:2 connected av:graphs}) implies \ref{thm:2 connected av:forcing}. The converse is taken care of by exhibiting certain constructions that show that \(321\), \(2341\) and \(4123\) fail to be \(2\)-forcing. We defer this to the appendix, Lemma~\ref{lem:not k-forcing}.

\begin{proof}[Proof of Theorem~\ref{thm:2 connected av}]
  We will prove that \(\ref{thm:2 connected av:forcing} \implies \ref{thm:2 connected av:patterns} \implies \ref{thm:2 connected av:graphs} \implies \ref{thm:2 connected av:forcing}\). First, suppose that \(\sigma\) contains one of the patterns \(\tau \in \{321, 2341, 4123\}\). Then \(\Av(\tau) \subseteq \Av(\sigma)\), so since condition \ref{thm:2 connected av:forcing} does not hold for \(\tau\) by Lemma~\ref{lem:not k-forcing} in the appendix, it also cannot hold for \(\sigma\). This proves the first implication.

  Second, suppose \(\sigma \in \Av(321, 2341, 4123)\). It is clear that \(G_{\sigma}\) contains no vertex of degree three or higher, so its components must be paths and cycles. \(3\)-cycles are impossible, since \(\sigma\) avoids \(321\). Cycles of length \(5\) or longer are ruled out by \(4\)-chordality.

  Lastly, decompose \(\sigma\) into its components \(\sigma = \sigma_1 \oplus \ldots \oplus \sigma_r\), and assume that each \(G_{\sigma_i}\) is a path or a \(4\)-cycle. By Lemmas~\ref{lem:3412 is 2-forcing} and \ref{lem:paths are 2-forcing}, each \(\sigma_i\) is \(2\)-forcing. By Theorem~\ref{thm:direct sum closure}, this implies that \(\sigma\) itself is \(2\)-forcing. The concrete bounds for \(C_2(\sigma_i)\) given by these results yield \(C_2(\sigma) \leq |\sigma| + 2p_1 + p_2 - q\), so we may take \(C = |\sigma| + 2p_1 + p_2 - q\) in \ref{thm:2 connected av:forcing}.
\end{proof}

\section[\texorpdfstring{To \(k\) with Menger}{To k with Menger}]{\texorpdfstring{To \(\bm k\) with Menger}{To k with Menger}} \label{sec:menger}

In this section we prove both Theorem~\ref{thm:k connected} and Theorem~\ref{thm:k connected av}. The revelation is that these theorems, stronger though they may appear, follow with relative ease from the \(k = 2\) results we already proved. We need only the well-known theorem of Menger \cite{menger_zur_1927} (see also Diestel \cite[Section~3.3]{diestel_graph_2025}, or any of the other standard textbooks) that there exist \(k\) pairwise disjoint paths between any two disjoint \(k\)-vertex sets in a \(k\)-connected graph. Throughout this section, we will assume the following setup.

\begin{construction} \label{const:menger}
  Let \(G_{\pi}\) be a \(k\)-connected inversion graph with \(n \geq 2k\) vertices. Denote by \(A = \{1, \ldots, k\}\) and \(B = \{n, \ldots, n-k+1\}\) the sets consisting of the bottom \(k\) and top \(k\) vertices, respectively. By Menger's theorem, there are \(k\) vertex-disjoint paths from \(A\) to \(B\). Let \(P_1, \ldots, P_k\) denote \(k\) such paths, each induced, ordered so that vertex \(i \in A\) is an endpoint of \(P_i\). Denote by \(H_{\pi}\) the subgraph induced on \(\bigcup_i P_i\). Write \(n_i\) for the number of vertices on \(P_i\), and \(m = \sum_i n_i\).
\end{construction}

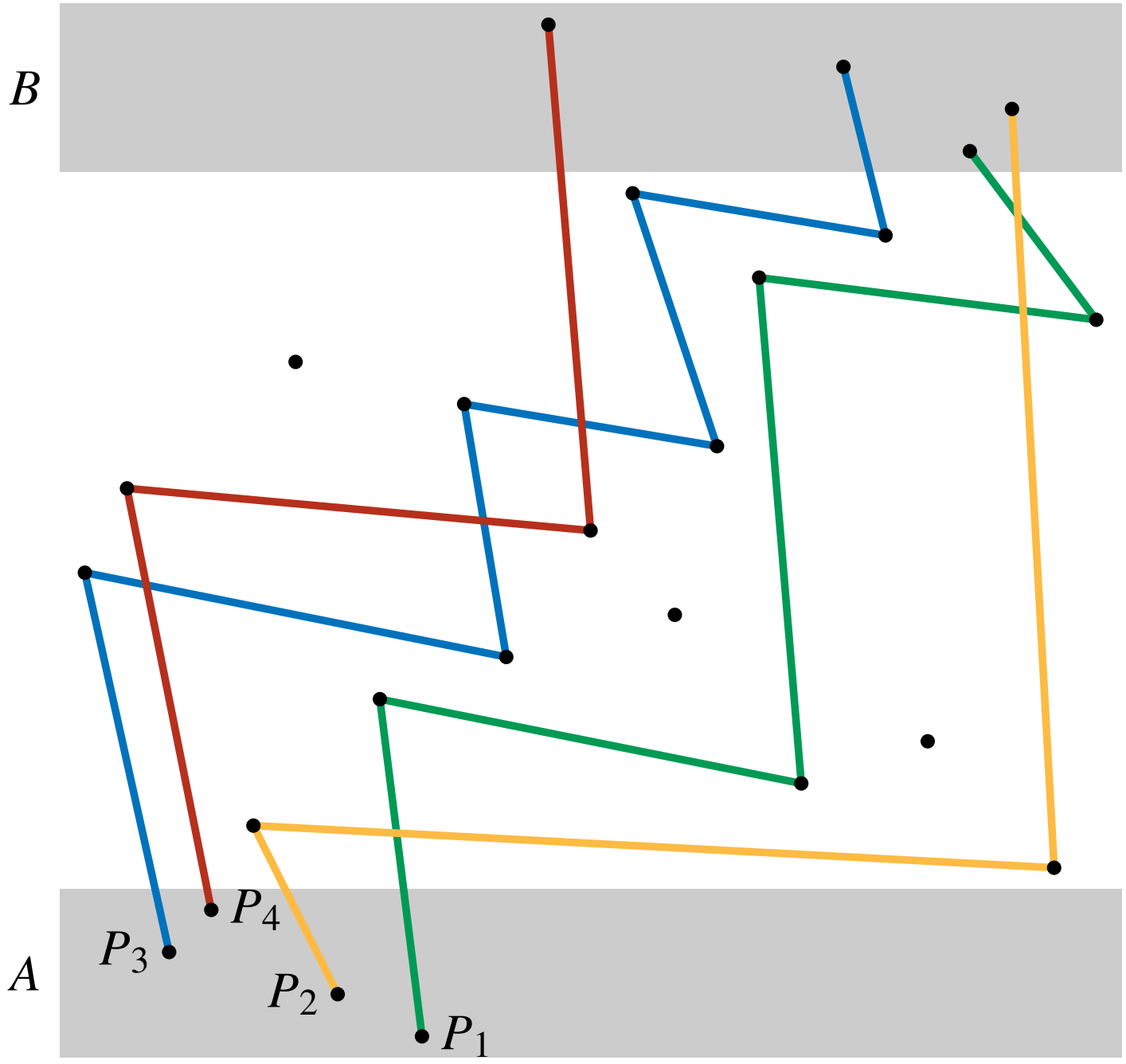
\begin{figure}[t]
  \centering
  \begin{tikzpicture}[scale=0.35]
    \fill[black!20] (0.4,0.5) rectangle (25.6,4.5);
    \fill[black!20] (0.4,21.5) rectangle (25.6,25.5);
    \node[left] at (0.4,2.5) {\(A\)};
    \node[left] at (0.4,23.5) {\(B\)};

    \foreach \v [count=\i] in {12,14,3,4,6,17,2,9,1,16,10,25,13,21,11,15,19,7,24,20,8,22,23,5,18}
      \coordinate (\v) at (\i,\v);

    \begin{scope}[line width=1.8pt, line cap=round, line join=round]
      \draw[ForestGreen] (1) -- (9) -- (7) -- (19) -- (18) -- (22);
      \draw[Dandelion] (2) -- (6) -- (5) -- (23);
      \draw[RoyalBlue] (3) -- (12) -- (10) -- (16) -- (15) -- (21) -- (20) -- (24);
      \draw[BrickRed] (4) -- (14) -- (13) -- (25);
    \end{scope}

    \foreach \v in {1,...,25}
      \node[dot] at (\v) {};

    \node[right] at (1) {\(P_1\)};
    \node[left] at (2) {\(P_2\)};
    \node[left] at (3) {\(P_3\)};
    \node[right] at (4) {\(P_4\)};
  \end{tikzpicture}
  \caption{Four disjoint induced paths \(P_1, \ldots, P_4\) from \(A\) to \(B\) in a \(4\)-connected inversion graph.}
  \label{fig:menger}
\end{figure}

Figure~\ref{fig:menger} gives an example of what this might look like when \(k = 4\). The paths \(P_1, \ldots, P_4\) are colored. Of course, there are many other edges in \(G_{\pi}\), which we have omitted for readability.

The structure of our argument is as follows:
\begin{itemize}
  \item We show that every vertex not in \(H_{\pi}\) has at least \(k\) edges to \(H_{\pi}\). Thus, it suffices to estimate edges inside \(H_{\pi}\).
  \item For Theorem~\ref{thm:k connected}, we can do this directly.
  \item For Theorem~\ref{thm:k connected av}, we need to estimate the number of edges between \(P_i\) and \(P_j\) for \(i \neq j\). We construct a \(2\)-connected inversion graph \(H_{ij}\) from \(P_i \cup P_j\). We use Theorem~\ref{thm:2 connected av} to estimate the number of edges inside \(H_{ij}\), and sum over all \(i < j\).
\end{itemize}

\begin{lemma} \label{lem:menger outside v}
  Every vertex in \(G_{\pi} \setminus (A \cup B)\) has at least one neighbor on each \(P_i\).
\end{lemma}

\begin{proof}
  Let \(v\) be a vertex in \(G_{\pi} \setminus (A \cup B)\), and assume that \(v \notin P_i\). Consider the induced subgraph \(G'\) on \(v \cup P_i\). The bottom and top vertices of \(G'\) are in the same component. Since \(G'\) is an inversion graph, this means that it must be connected. In particular, \(v\) must have a neighbor on \(P_i\).
\end{proof}

This is all we need for Theorem~\ref{thm:k connected}. Let us recall also the statement we wish to prove.

\begin{retheorem}{thm:k connected}
  Fix \(k \geq 1\). There exists a constant \(C\) depending only on \(k\) such that every permutation \(\pi\) whose inversion graph is \(k\)-connected satisfies
  \begin{equation*}
    \inv(\pi) \geq \frac{k+1}{2} |\pi| - C.
  \end{equation*}
\end{retheorem}

\begin{proof}
  Let \(\pi \in S_n\) be a permutation whose inversion graph is \(k\)-connected. First, assume \(n \geq 2k\) and perform Construction~\ref{const:menger}.
  
  By Lemma~\ref{lem:menger outside v}, the internal vertices of the path \(P_i\) all have a neighbor in every other \(P_j\). Since they also have two neighbors each inside their own path, the degree of such a vertex in \(H_{\pi}\) is at least \(k+1\). 

  On the other hand, the bottom endpoint \(b\) of \(P_i\) must have a neighbor on each path \(P_j\) whose bottom endpoint is to the left of \(b\), since \(P_j\) must cross the region northwest of \(b\). In other words, the vertices of \(A\) have degrees at least \(1, \ldots, k\) in \(H_{\pi}\), and the same is true for \(B\). Summing degrees, we find that
  \begin{equation*}
    2|E(H_{\pi})| \geq (k+1)(m - 2k) + 2 \sum_{i = 1}^k i = (k+1)m - k(k+1),
  \end{equation*}
  and thus
  \begin{equation*}
    |E(H_{\pi})| \geq \frac{k+1}{2} m - \binom{k+1}{2}.
  \end{equation*}
  A second application of Lemma~\ref{lem:menger outside v} shows that every vertex in \(G_{\pi} \setminus H_{\pi}\) has at least one neighbor on each \(P_i\). We get
  \begin{equation*}
    |E(G_{\pi})| \geq |E(H_{\pi})| + k(n - m) \geq \frac{k+1}{2}n - \binom{k + 1}{2}.
  \end{equation*}
  Hence, Theorem~\ref{thm:k connected} holds with constant \(C = \binom{k+1}{2}\).

  On the other hand, if \(n \leq 2k\), we find that even with the smaller constant \(C = k\),
  \begin{equation*}
    \frac{k+1}{2}n - k \leq \frac{k+1}{2}n - \frac n2 = \frac{kn}{2} \leq |E(G_{\pi})|.
  \end{equation*}
  In either case, the bound holds.
\end{proof}

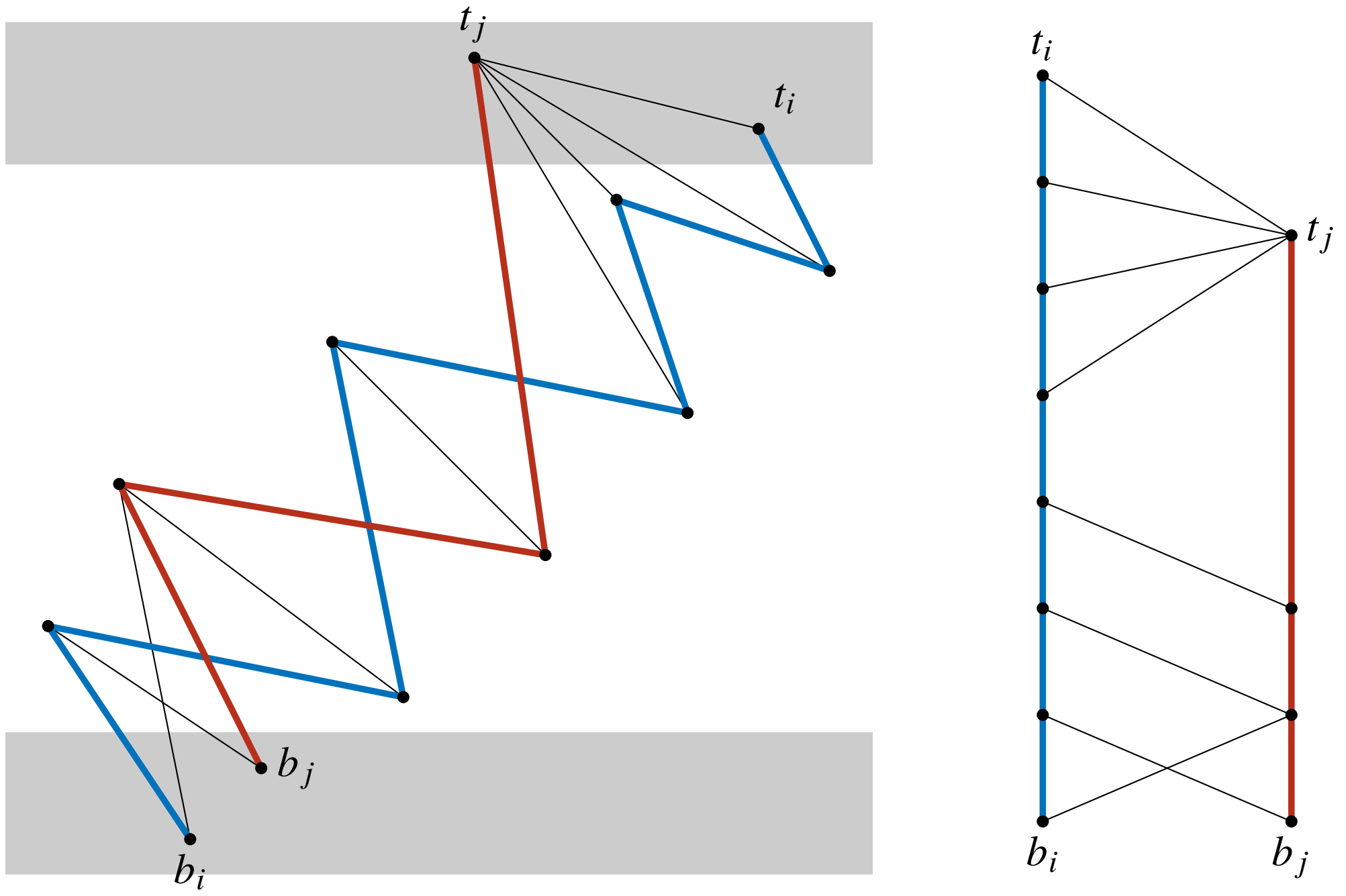
\begin{figure}[t]
  \centering
  \begin{tikzpicture}[scale=0.7]
    \def\twopaths{
      \draw[rpath]
        (1) -- (6)  (2) -- (4)  (3) -- (6)  (5) -- (8)
        (7) -- (12)  (9) -- (12)  (10) -- (12)  (11) -- (12);
      \begin{scope}[line width=1.8pt, line cap=round, line join=round]
        \draw[RoyalBlue] (1) -- (4) -- (3) -- (8) -- (7) -- (10) -- (9) -- (11);
        \draw[BrickRed] (2) -- (6) -- (5) -- (12);
      \end{scope}
      \foreach \v in {1,...,12}
        \node[dot] at (\v) {};
    }

    \fill[black!20] (0.4,0.5) rectangle (12.6,2.5);
    \fill[black!20] (0.4,10.5) rectangle (12.6,12.5);
    \foreach \v [count=\i] in {4,6,1,2,8,3,12,5,10,7,11,9}
      \coordinate (\v) at (\i,\v);
    \twopaths
    \node[below] at (1) {\(b_i\)};
    \node[right] at (2) {\(b_j\)};
    \node[above right] at (11) {\(t_i\)};
    \node[above] at (12) {\(t_j\)};

    \begin{scope}[shift={(15,2.75)}]
      \foreach \v/\x/\y in {1/0/-1.5, 4/0/0, 3/0/1.5, 8/0/3, 7/0/4.5, 10/0/6, 9/0/7.5, 11/0/9,
          2/3.5/-1.5, 6/3.5/0, 5/3.5/1.5, 12/3.5/6.75}
        \coordinate (\v) at (\x,\y);
      \twopaths
      \node[below] at (1) {\(b_i\)};
      \node[below] at (2) {\(b_j\)};
      \node[above] at (11) {\(t_i\)};
      \node[right] at (12) {\(t_j\)};
    \end{scope}
  \end{tikzpicture}
  \caption{The subgraph induced by the two leftmost paths \(P_3\) and \(P_4\) of Figure~\ref{fig:menger}: a diagram representing it (left), and its graph-theoretical structure (right).}
  \label{fig:menger two paths}
\end{figure}

Let us now fix some \(i \neq j\). The subgraph induced on \(P_i \cup P_j\) consists of two paths with some edges in between. Denote the bottom and top endpoints of \(P_i\) by \(b_i\) and \(t_i\), respectively. We define \(H_{ij}\) as the graph obtained from \(G_{\pi}[P_i \cup P_j]\) by including the new edges \(b_i b_j\) and \(t_i t_j\), unless these edges already exist. Figure~\ref{fig:menger two paths} shows the induced subgraph formed by the paths \(P_3\) and \(P_4\) from Figure~\ref{fig:menger}. To form \(H_{ij}\), we include an edge connecting the bottom two vertices. This corresponds to swapping the vertical coordinates of the bottom two points in the diagram on the left.

\begin{lemma} \label{lem:menger two paths 2 connected}
  The graph \(H_{ij}\) is a \(2\)-connected inversion graph. 
\end{lemma}

\begin{proof}
  \(2\)-connectivity is clear: since the endpoints of the two paths are connected, \(H_{ij}\) contains a Hamiltonian cycle. To see that \(H_{ij}\) is an inversion graph, observe that a diagram for \(H_{ij}\) can be obtained from a diagram of \(G_{\pi}[P_i \cup P_j]\) by possibly swapping the vertical coordinates of \(b_i\) with \(b_j\) and \(t_i\) with \(t_j\) so that both pairs form inversions.
\end{proof}
 
This lemma means that we can apply our \(k=2\) results to \(H_{ij}\). However, we must take care to use the right constants, since shifting around the endpoints might have created some new patterns.

We are ready for the interesting part of Theorem~\ref{thm:k connected av}. The easy direction, \ref{rethm:k connected av:forcing} implies \ref{rethm:k connected av:patterns}, is taken care of by Lemma~\ref{lem:not k-forcing} in the appendix. Let us restate and prove our main result.

\begin{retheorem}{thm:k connected av}
  Fix \(k \geq 2\), and let \(\sigma\) be a pattern. The following are equivalent:
  \begin{enumerate}[(i)]
    \item\label{rethm:k connected av:forcing} There exists a constant \(C\) depending only on \(k\) and \(\sigma\) such that every permutation \(\pi \in \Av(\sigma)\) whose inversion graph is \(k\)-connected satisfies
    \begin{equation*}
      \inv(\pi) \geq k|\pi| - C.
    \end{equation*}
    \item\label{rethm:k connected av:patterns} We have \(\sigma \in \Av(321, 2341, 4123)\).
    \item\label{rethm:k connected av:graphs} The components of \(G_{\sigma}\) are paths and \(4\)-cycles.
  \end{enumerate}
\end{retheorem}

\begin{proof}
  Theorem~\ref{thm:2 connected av} establishes \(\ref{rethm:k connected av:patterns} \iff \ref{rethm:k connected av:graphs}\), and Lemma~\ref{lem:not k-forcing} gives \(\ref{rethm:k connected av:forcing} \implies \ref{rethm:k connected av:patterns}\), so only one implication remains. Let \(\sigma \in \Av(321, 2341, 4123)\), and assume that \(\pi \in \Av_n(\sigma)\) has a \(k\)-connected inversion graph. As before, assume that \(n \geq 2k\), and perform Construction~\ref{const:menger}.

  Fix \(i \neq j\). Our graph \(H_{ij}\) might contain \(\sigma\), but we claim that it cannot contain the pattern \(\tau = 1 \oplus \sigma \oplus 1\). Indeed, \(H_{ij}\) was built from the subgraph induced on \(P_i \cup P_j\) by (possibly) swapping \(b_i\) with \(b_j\) and \(t_i\) with \(t_j\), so that both pairs form inversions. The two pairs of points are, respectively, the bottom two points and top two points in \(H_{ij}\). It follows that if an occurrence of \(\tau\) in \(H_{ij}\) uses at most one entry from each of the sets \(\{b_i, b_j\}\) and \(\{t_i, t_j\}\), then the same pattern can be found in \(G_{\pi}\). Thus, without loss of generality, we may assume that the occurrence contains both \(b_i\) and \(b_j\). However, these two points must form an inversion, whereas the entries \(1\) and \(2\) of \(\tau\) do not. This is a contradiction, so \(H_{ij}\) must avoid \(\tau\).

  Since \(C_2(1) = 3\) by Lemma~\ref{lem:paths are 2-forcing}, Theorems~\ref{thm:direct sum closure} and \ref{thm:2 connected av} show that our new pattern \(\tau\) is \(2\)-forcing with constant
  \begin{equation*}
    C_2(\tau) \leq C_2(\sigma) + 6 \leq |\sigma| + 2p_1 + p_2 - q + 6,
  \end{equation*}
  where \(p_1\), \(p_2\) and \(q\) count the number of isolated vertices, isolated edges and \(4\)-cycles of \(G_{\sigma}\), respectively. Thus, with Lemma~\ref{lem:menger two paths 2 connected},
  \begin{equation*}
    |E(H_{ij})| \geq 2(n_i + n_j) - C_2(\sigma) - 6.
  \end{equation*}
  Writing \(e_{ij}\) for the number of edges between \(P_i\) and \(P_j\) in \(G_{\pi}\), we have
  \begin{equation*}
    (n_i - 1) + (n_j - 1) + e_{ij} + 2 \geq |E(H_{ij})|,
  \end{equation*}
  which, together with the previous inequality, yields
  \begin{equation*}
    e_{ij} \geq n_i + n_j - C_2(\sigma) - 6.
  \end{equation*}

  Summing over all pairs \(i < j\) and including the edges of each \(P_i\) gives
  \begin{align*}
    |E(H_{\pi})| &\geq m - k + \sum_{i < j} e_{ij} \\
    &\geq m - k + (k - 1)m - \binom k2 (C_2(\sigma) + 6) \\
    &= km - k - \binom k2 (C_2(\sigma) + 6).
  \end{align*}
  As before, this yields
  \begin{equation*}
    |E(G_{\pi})| \geq kn - k - \binom k2 (C_2(\sigma) + 6), \postdisplaypenalty=10000
  \end{equation*}
  proving Theorem~\ref{thm:k connected av} with constant \(C = k + \binom k2 (C_2(\sigma) + 6)\).

  If, on the other hand, \(n \leq 2k\), then we can take the smaller constant \(C = k^2\), since
  \begin{equation*}
    kn - k^2 \leq kn - \frac{kn}{2} = \frac{kn}{2} \leq |E(G_{\pi})|.
  \end{equation*}
  In both cases, the desired bound stands.
\end{proof}

\section[\texorpdfstring{Weak \(2\)-connectivity}{Weak 2-connectivity}]{\texorpdfstring{Weak \(\bm 2\)-connectivity}{Weak 2-connectivity}} \label{sec:weak 2-connectivity}

In this section, we return to the original motivation for this paper: the result of Linusson and Verkama \cite{linusson_enumerating_2025} that every \(1324\)-avoiding permutation of size \(n\) with at most \(2n - 7\) inversions must be either decomposable or \emph{almost} decomposable. We will interpret this in terms of the inversion graph.

\begin{definition}
  An inversion graph \(G_{\pi}\) on \([n]\) is called \emph{weakly \(2\)-connected} if \(G_{\pi}\) is connected and
  \begin{equation} \label{eq:weakly 2-connected}
    G_{\pi} \setminus v \text{ is connected for every } v \in \{1, n, \pi(1), \pi(n)\}.
  \end{equation}
\end{definition}

\begin{remark}
  Be careful: this is not an abstract graph-theoretical property, since condition \eqref{eq:weakly 2-connected} specifies a labeled set of vertices. In particular, it is possible that two inversion graphs \(G_{\pi}\) and \(G_{\tau}\) are isomorphic, but only one of the two is weakly \(2\)-connected. 
\end{remark}

\begin{remark}
  An inversion graph \(G_{\pi}\) is weakly \(2\)-connected if and only if \(\pi\) is neither decomposable nor almost decomposable. Thus, every \(1324\)-avoiding weakly \(2\)-connected inversion graph has at least \(2n - 6\) edges by \cite[Theorem~8]{linusson_enumerating_2025}.
\end{remark}

\begin{definition}
  A pattern \(\sigma\) is called \emph{strongly \(2\)-forcing} if there is a constant \(C\) depending only on \(\sigma\) such that every permutation \(\pi \in \Av(\sigma)\) whose inversion graph is weakly \(2\)-connected satisfies
  \begin{equation*}
    \inv(\pi) \geq 2|\pi| - C.
  \end{equation*}
\end{definition}

\begin{remark}
  Every strongly \(2\)-forcing pattern is \(2\)-forcing by definition. Hence the name.
\end{remark}

Our proof of Theorem~\ref{thm:weakly 2-connected} revolves around an idea developed in a recent AI-assisted preprint of Allikvere \cite{allikvere_inversion_2026}; see particularly their Section~6.2. The main observation is that the diameter of an inversion graph avoiding a path is bounded, and this bounds the size of certain induced subgraphs.

\begin{definition}
  An induced subgraph \(H\) of an inversion graph \(G_{\pi}\) is called \emph{boundary-spanning} if \(\{1, n, \pi(1), \pi(n)\} \subseteq V(H)\).
\end{definition}

We will construct a specific such subgraph. Note the similarity of the following with Construction~\ref{const:menger} from the previous section. In this part, we will (mostly for convenience) assume that all four entries \(1\), \(n\), \(\pi(1)\) and \(\pi(n)\) are distinct. The other cases are not difficult, and we will handle them later.

\begin{construction} \label{const:boundary pattern}
  Let \(G_{\pi}\) be a weakly \(2\)-connected inversion graph on \([n]\) such that all four entries \(1\), \(n\), \(\pi(1)\) and \(\pi(n)\) are distinct. Let \(P_1\) and \(P_2\) be induced paths from \(\pi(1)\) to \(\pi(n)\) in \(G_{\pi} \setminus 1\) and \(G_{\pi} \setminus n\), respectively. Similarly, let \(Q_1\) and \(Q_2\) be induced paths from \(1\) to \(n\) in \(G_{\pi} \setminus \pi(1)\) and \(G_{\pi} \setminus \pi(n)\), respectively. Define \(B_{\pi}\) as the subgraph of \(G_{\pi}\) induced on \(P_1 \cup P_2 \cup Q_1 \cup Q_2\).
\end{construction}

A critical part of the argument is that each strongly \(2\)-forcing pattern \(\sigma\) in Theorem~\ref{thm:weakly 2-connected} is contained in a path. This provides us with some unifying structure.

\begin{lemma}[Cf.\ Lemmas 6.2 and 6.3 in \cite{allikvere_inversion_2026}] \label{lem:boundary spanning}
  Every induced subgraph \(H\) of \(G_{\pi}\) containing \(B_{\pi}\) is weakly \(2\)-connected. Furthermore, if \(\pi\) avoids a pattern \(\sigma\) whose inversion graph is a path on \(m \geq 3\) vertices, then \(|V(B_{\pi})| \leq 4m - 6\).
\end{lemma}

\begin{proof}
  Observe that by construction, \(P_1\) is still contained in the inversion graph \(H \setminus 1\). It is a path connecting the leftmost and rightmost vertices of the inversion graph \(H\), so \(H \setminus 1\) must be connected. The other cases are symmetrical, so \(H\) is weakly \(2\)-connected.

  Suppose now that \(\pi\) avoids a pattern \(\sigma\) whose inversion graph is a path on \(m \geq 3\) vertices. Recall from before (and from \cite{albert_rationality_2019}) that there are exactly two such patterns, distinguished by whether the path begins at \(1\) or \(\sigma(1)\). Furthermore, both of the patterns contain all strictly shorter patterns whose inversion graphs are paths. It follows that two of our four paths (either \(P_1\) and \(P_2\) or \(Q_1\) and \(Q_2\)) have at most \(m - 1\) vertices, and the other two at most \(m\). Thus,
  \begin{equation*}
    |V(B_{\pi})| \leq 4 + 2(m - 3) + 2(m - 2) = 4m - 6. \qedhere
  \end{equation*}
\end{proof}

The idea for the rest of the proof is to show that under our pattern-avoidance conditions, the vertices in \(G_{\pi} \setminus B_{\pi}\) must have at least two neighbors in \(B_{\pi}\) (again, this is similar to part of the strategy in the previous section). This is not always literally true, but it is close enough.

From now on, we will assume that \(\pi\) is such that
\begin{equation} \label{eq:boundary condition}
  \pi(1) < \pi(n) \quad \text{and} \quad \pi^{-1}(1) < \pi^{-1}(n).
\end{equation}
(Observe that together with connectivity, this implies that all four corresponding entries are distinct, as we assumed in Construction~\ref{const:boundary pattern}.) This is the critical case, as Lemma~\ref{lem:trivial boundary condition} below from \cite{verkama_inversions_2026} shows that if it fails, then the permutation must have enough inversions, even without any pattern avoidance. Also, we will see that a connected \(2413\)-avoiding permutation cannot satisfy \eqref{eq:boundary condition}, so that case is taken care of automatically.

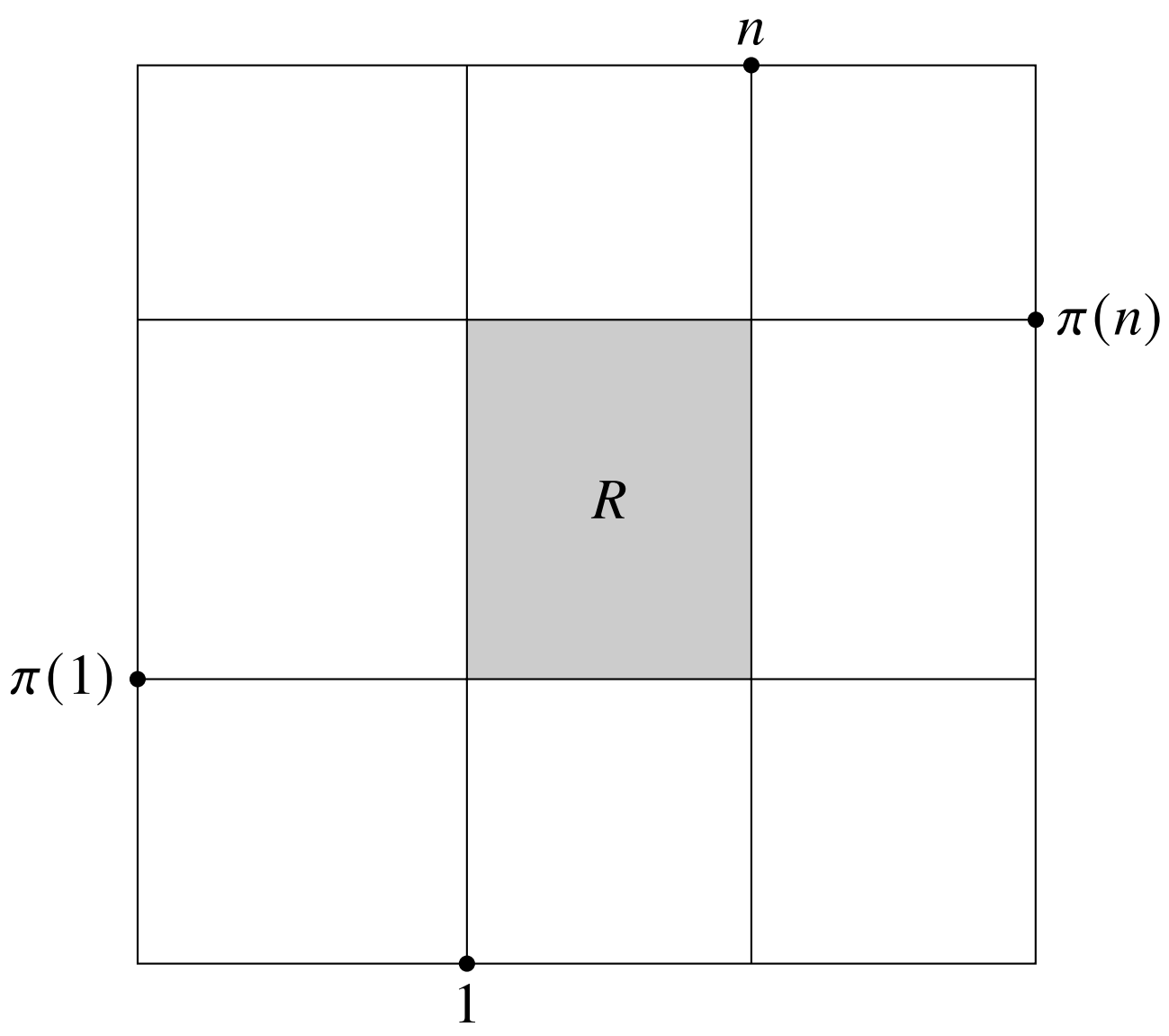
\begin{figure}[t]
  \centering
  \begin{tikzpicture}[scale=1.1]
    \coordinate (b) at (2.2,0);
    \coordinate (t) at (4.1,6);
    \coordinate (l) at (0,1.9);
    \coordinate (r) at (6,4.3);

    \fill[black!20] (b |- l) rectangle (t |- r);
    \draw (0,0) rectangle (6,6);
    \draw[rpath]
      (b) -- (b |- t)
      (t |- b) -- (t)
      (l) -- (r |- l)
      (l |- r) -- (r);

    \node at (3.15,3.1) {\(R\)};
    \foreach \l in {b, t, l, r}
      \node[dot] at (\l) {};
    \node[below] at (b) {\(1\)};
    \node[above] at (t) {\(n\)};
    \node[left] at (l) {\(\pi(1)\)};
    \node[right] at (r) {\(\pi(n)\)};
  \end{tikzpicture}
  \caption{The central region \(R\) in the diagram of \(G_{\pi}\), when \(\pi(1) < \pi(n)\) and \(\pi^{-1}(1) < \pi^{-1}(n)\).}
  \label{fig:central region}
\end{figure}

Assumption \eqref{eq:boundary condition} means that the diagram of \(G_{\pi}\) has a central region \(R = (\pi^{-1}(1), \pi^{-1}(n)) \times (\pi(1), \pi(n))\); see Figure~\ref{fig:central region}. The following lemmas develop our strategy.

\begin{lemma} \label{lem:v outside R has two neighbors}
  Let \(v \in G_{\pi} \setminus B_{\pi}\) be a vertex that is not in \(R\). Then \(v\) must have at least two neighbors in \(B_{\pi}\).
\end{lemma}

\begin{proof}
  It is easy to see by inspecting Figure~\ref{fig:central region} that since \(v\) is not in \(R\), it must have a neighbor \(u \in \{1, n, \pi(1), \pi(n)\}\). By Lemma~\ref{lem:boundary spanning}, the induced subgraph \(H = G_{\pi}[B_{\pi} \cup v]\) is weakly \(2\)-connected. Hence, \(H \setminus u\) is connected, and \(v\) must have some other neighbor in \(B_{\pi}\).
\end{proof}

\begin{lemma} \label{lem:R is empty}
  If \(\pi\) avoids \(1324\) or \(1342\), then \(R\) must be empty.
\end{lemma}

\begin{proof}
  Suppose \(R\) is nonempty. Since \(G_{\pi}\) is connected, some vertex \(v\) in \(R\) must have a neighbor \(u\) outside of \(R\). If this neighbor \(u\) is northwest of \(v\), then \(\pi(1) u v n\) is a \(1324\)-pattern. If \(u\) is southeast of \(v\), then \(1 v u \pi(n)\) is a \(1324\)-pattern. 

  It remains to show that \(\pi\) must contain \(1342\). Consider the diagram of \(G_{\pi}\) in Figure~\ref{fig:1342 central region}. First, if there is a point \(p\) in the upper grey region, then \(\pi(1) p n \pi(n)\) is a \(1342\)-pattern. Similarly, a point \(q\) in the lower grey region creates the \(1342\)-pattern \(1 v n q\). Assume that both are empty.

  \begin{figure}[t]
    \centering
    \begin{tikzpicture}[scale=1.1]
      \coordinate (b) at (2.2,0);
      \coordinate (t) at (4.1,6);
      \coordinate (l) at (0,1.9);
      \coordinate (r) at (6,4.3);
      \coordinate (v) at (3.15,3.1);
      \coordinate (x) at (1.1,3.7);
      \coordinate (y) at (5.05,3.45);

      \fill[black!20]
        (l |- r) rectangle (t)
        (t |- b) rectangle (r |- l)
        (t |- l) rectangle (r |- v);
      \fill[RoyalBlue!30] (l |- v) rectangle (t |- r);
      \fill[BrickRed!30] (t |- v) rectangle (r);
      \draw (0,0) rectangle (6,6);
      \draw[rpath]
        (b) -- (b |- t)
        (t |- b) -- (t)
        (l) -- (r |- l)
        (l |- r) -- (r)
        (l |- v) -- (r |- v);
      \draw[rpath, thick] (x) -- (y);

      \foreach \l in {b, t, l, r, x, v, y}
        \node[dot] at (\l) {};
      \node[below] at (b) {\(1\)};
      \node[above] at (t) {\(n\)};
      \node[left] at (l) {\(\pi(1)\)};
      \node[right] at (r) {\(\pi(n)\)};
      \node[left] at (x) {\(x\)};
      \node[below] at (v) {\(v\)};
      \node[right] at (y) {\(y\)};
    \end{tikzpicture}
    \caption{The situation in the proof of Lemma~\ref{lem:R is empty}.}
    \label{fig:1342 central region}
  \end{figure}
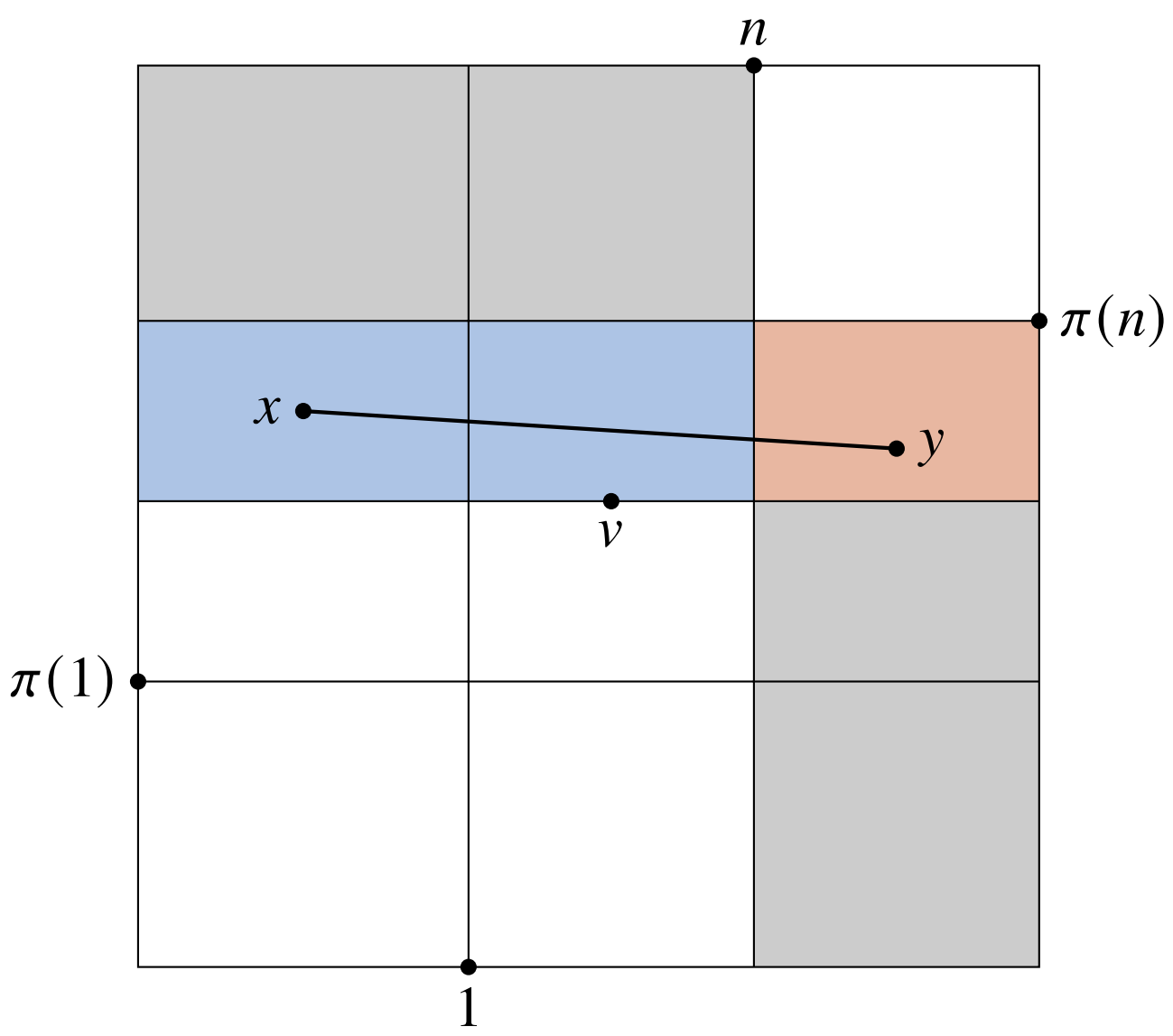

  Now, it is clear that for \(G_{\pi}\) to be connected, it must contain an edge \(xy\) from the blue region to the red. Then \(\pi(1) x n y\) is a \(1342\)-pattern, concluding the proof.
\end{proof}

The last patterns to consider, as per the theorem statement, are those of the form \(21 \oplus {\id_t} \oplus 21\). These are not quite as simple, as \(R\) may not actually be empty. However, it is clearly true that \(R\) cannot contain an increasing sequence of size \(t\). This is sufficient.

\begin{lemma} \label{lem:no long increase in R}
  If \(\pi\) avoids \(21 \oplus {\id_t} \oplus 21\), then \(R\) cannot contain an increasing sequence of size \(t\).
\end{lemma}

We are now almost ready for the hard direction of Theorem~\ref{thm:weakly 2-connected}. The following lemma takes care of when the condition \eqref{eq:boundary condition} fails. This is \cite[Lemma~C.1.6]{verkama_inversions_2026}, a simple generalization of \cite[Lemma~14]{linusson_enumerating_2025}.

\begin{lemma}[Lemma C.1.6 in \cite{verkama_inversions_2026}] \label{lem:trivial boundary condition}
  If \(G_{\pi}\) is a weakly \(2\)-connected inversion graph on \([n]\) and \(\pi(1) > \pi(n)\) or \(\pi^{-1}(1) > \pi^{-1}(n)\), then \(\inv(\pi) \geq 2n - 4\). 
\end{lemma}

For the other implication in Theorem~\ref{thm:weakly 2-connected}, we need to show that no other patterns are strongly \(2\)-forcing. This is done by identifying an explicit collection \(F\) of patterns such that the permutation class \(\Sub(\Gamma B)\) appearing in the theorem statement equals \(\Av(F)\). We do this in Lemma~\ref{lem:weakly 2-forcing class} of the appendix. We are now prepared.

\begin{retheorem}{thm:weakly 2-connected}
  Let \(\sigma\) be a pattern. The following are equivalent:
  \begin{enumerate}[(i)]
    \item\label{rethm:weakly 2-connected:forcing} There exists a constant \(C\) depending only on \(\sigma\) such that every \(\pi \in \Av(\sigma)\) whose inversion graph is weakly \(2\)-connected satisfies
    \begin{equation*}
      \inv(\pi) \geq 2|\pi| - C.
    \end{equation*}
    \item\label{rethm:weakly 2-connected:patterns} We have \(\sigma \in \Sub(\Gamma B)\), where
    \begin{equation*}
      B = \left\{1324, 1342, 2413\right\} \cup \left\{21 \oplus {\id_t} \oplus 21 : t \in \mathbb Z_{\geq 0}\right\}.
    \end{equation*}
  \end{enumerate}
\end{retheorem}

\begin{proof}
  The implication \(\ref{rethm:weakly 2-connected:forcing} \implies \ref{rethm:weakly 2-connected:patterns}\) follows directly from Lemma~\ref{lem:weakly 2-forcing class}, so it suffices to prove the converse. Since the set of strongly \(2\)-forcing patterns is closed under pattern containment as well as the action of \(\Gamma\), it is enough to prove that every pattern in \(B\) is strongly \(2\)-forcing.

  Let \(\pi \in S_n\) be a permutation whose inversion graph is weakly \(2\)-connected. If we have that \(\pi(1) > \pi(n)\) or \(\pi^{-1}(1) > \pi^{-1}(n)\), then Lemma~\ref{lem:trivial boundary condition} proves that \(\inv(\pi) \geq 2n - 4\), and we are done even without any pattern-avoidance constraints. So, suppose that \(\pi(1) < \pi(n)\) and \(\pi^{-1}(1) < \pi^{-1}(n)\), and carry out Construction~\ref{const:boundary pattern}.

  First, let \(\sigma = 2413\), and let \(P\) be an induced path from \(\pi(1)\) to \(\pi(n)\). Since \(P\) begins at the first entry and ends at the last entry, it must contain an even number of vertices. Since \(\pi(1) < \pi(n)\), this number is at least four. However, then \(P\) contains \(2413\), so \(\pi\) cannot avoid \(2413\). This proves that \(\sigma = 2413\) satisfies condition \ref{rethm:weakly 2-connected:forcing} with the constant \(C = 4\) from Lemma~\ref{lem:trivial boundary condition}.

  Take instead \(\sigma \in \{1324, 1342\}\), and suppose that \(\pi \in \Av(\sigma)\). By Lemmas~\ref{lem:v outside R has two neighbors} and \ref{lem:R is empty}, every vertex of \(G_{\pi} \setminus B_{\pi}\) must have at least two neighbors in \(B_{\pi}\). Furthermore, \(\sigma\) is contained in the pattern \(241635\), whose inversion graph is a path on \(6\) vertices. By Lemma~\ref{lem:boundary spanning}, we must have \(|V(B_{\pi})| \leq 18\). All taken together, this yields
  \begin{equation*}
    \inv(\pi) \geq 2(n - 18) + |E(B_{\pi})| \geq 2n - 36.
  \end{equation*}
  Hence, the claim holds with \(C = 36\) (obviously, this is far from optimal).

  Our last case is \(\sigma = 21 \oplus {\id_t} \oplus 21\). If \(t = 0\), then \(\sigma = 2143\), and \(\pi\) contains it: \(\pi(1) 1 n \pi(n)\) is a \(2143\)-pattern. We may therefore assume that \(t \geq 1\). This pattern \(\sigma\) is contained in both permutations of length \(2t + 6\) whose inversion graphs are paths, so \(|V(B_{\pi})| \leq 8t + 18\). Now, denote by \(R'\) the set of vertices of \(G_{\pi} \setminus B_{\pi}\) that are contained in the central region \(R\) (see Figure~\ref{fig:central region}).

  Let us construct an induced subgraph \(H\) of \(G_{\pi}\) by inserting the points of \(R'\) into \(B_{\pi}\) one by one. By Lemma~\ref{lem:no long increase in R}, \(R'\) cannot contain an increasing sequence of size \(t\). In particular, it can happen at most \(t-1\) times that a point we insert has no neighbors among the previously inserted points of \(R'\). Since \(B_{\pi}\) is connected and boundary-spanning, each vertex also has a neighbor in \(B_{\pi}\). This means that all but \(t-1\) of the points of \(R'\) create at least two new edges. For the points outside \(R'\), we can, again, use Lemma~\ref{lem:v outside R has two neighbors}. All in all,
  \begin{align*}
    \inv(\pi) &\geq |E(B_{\pi})| + 2(n - |V(B_{\pi})|) - (t - 1) \\
    &\geq 2(n - 8t - 18) - t + 1 = 2n - 17t - 35.
  \end{align*}
  The statement holds with \(C = 17t + 35\).
\end{proof}

\section{Future directions} \label{sec:future work}

\paragraph{On tightness}

It is not too hard to construct examples showing that for every \(k \geq 2\), the linear coefficient \(k\) in the bound
\begin{equation*}
  \inv(\pi) \geq k|\pi| - C
\end{equation*}
in Theorem~\ref{thm:k connected av} \ref{thm:k connected av:forcing} is generally tight. The same goes for the corresponding coefficient \(2\) in Theorem~\ref{thm:weakly 2-connected}. However, we did not do too much work to optimize the constants \(C\). Concretely, our proof of Theorem~\ref{thm:weakly 2-connected} yields the constant \(C = 36\) for the pattern \(1324\), whereas the best possible constant is known to be \(6\) \cite[Theorem~8]{linusson_enumerating_2025}.

Theorems~\ref{thm:k connected} and \ref{thm:k connected av} are even worse. In Theorem~\ref{thm:k connected}, our proof yields \(C = \binom{k+1}{2}\), whereas data suggest that \(C = k\) should work. The other result relies on the same method, so the constant we obtained in its proof should also be far from optimal.

The \(2\)-connected versions are better. The bound of Proposition~\ref{prop:l chordal bound} for \(\ell\)-chordal graphs is tight (up to rounding) for every \(\ell\), and so is its restriction to inversion graphs (a subclass of the \(\ell = 4\) case). We also expect the bound given by Theorem~\ref{thm:direct sum closure} for direct sums of patterns to be close to optimal. Thus, the constant \(C\) stated in Theorem~\ref{thm:2 connected av} should not be too far from the truth.

\begin{problem}
  Optimize the constants \(C\) in Theorems~\ref{thm:k connected}, \ref{thm:k connected av} and \ref{thm:weakly 2-connected}.
\end{problem}

\paragraph{Weaker forcing conditions}

Since the ladder graphs are the unique \(2\)-connected inversion graphs on \(2n\) vertices with the minimum possible number of edges, it makes sense that they need to be avoided to increase the bound of Theorem~\ref{thm:k connected} (for \(k = 2\)); and indeed, our patterns \(\sigma \in \Av(321, 2341, 4123)\) are contained in the ladder graph.

However, it is not so obvious what happens if different pieces of the ladder are avoided. Remark~\ref{rmk:weaker forcing example} in the appendix discusses a certain family of permutations \(\tau\) with
\begin{equation*}
  \inv(\tau) = \left\lceil\frac 74 |\tau|\right\rceil - 3.
\end{equation*}
An example is visualized in Figure~\ref{fig:tau big}. The permutation \(\tau\) avoids \(2341\) and thus also \(351624\), whose inversion graph is the \(3\)-rung ladder graph. Hence, for \(\sigma = 351624\), the best result we can hope for is that every \(\pi \in \Av(\sigma)\) whose inversion graph is \(2\)-connected satisfies \(\inv(\pi) \geq \frac 74 |\pi| - 3\).
If this is true, it would be evidence of some richer spectrum of structure than what we found in this article. Our computations suggest that it is true.

\begin{conjecture}
  Every permutation \(\pi \in \Av(351624)\) whose inversion graph is \(2\)-connected satisfies
  \begin{equation*}
    \inv(\pi) \geq \frac 74 |\pi| - 3.
  \end{equation*}
\end{conjecture}

One could imagine that this generalizes to even longer ladders in some way.

\paragraph{Connectivity as a statistic}

Denote by \(\kappa(G)\) the largest \(k\) such that the graph \(G\) is \(k\)-connected. In light of our results, it makes sense to study the distribution of \(\kappa\) over inversion graphs. Define
\begin{equation*}
  K_n(x) = \sum_{\pi \in S_n} x^{\kappa(G_{\pi})}.
\end{equation*}
The triangle of coefficients of \(K_n(x)\) for \(n \leq 7\) is given in Table~\ref{tab:connectivity}. It does not appear to be in the OEIS \cite{oeis}. Some parts of it are familiar: the first column (for \(n \geq 2\)), i.e.\ the constant term of \(K_n(x)\), is the number of decomposable permutations, sequence \href{https://oeis.org/A356291}{A356291}. The last term of each row is \(1\), corresponding to the decreasing permutation (i.e.\ the complete graph). The polynomials \(K_n(x)\) are not real-rooted in general: this fails already for \(n = 3\). Log-concavity and unimodality hold in the rows we computed.

\begin{table}[t]
  \centering
  \caption{The number of permutations \(\pi \in S_n\) with \(\kappa(G_{\pi}) = k\), for \(n \leq 7\) (rows) and \(0 \leq k \leq n - 1\) (columns).}
  \label{tab:connectivity}
  \begin{tabular}{rrrrrrr}
    1 \\
    1 & 1 \\
    3 & 2 & 1 \\
    11 & 8 & 4 & 1 \\
    49 & 42 & 21 & 7 & 1 \\
    259 & 256 & 137 & 55 & 12 & 1 \\
    1593 & 1778 & 1065 & 448 & 135 & 20 & 1
  \end{tabular}
\end{table}

It is noteworthy that the distribution is so right-skewed. Bollobás and Thomason \cite{bollobas_random_1985} showed that in the Erd\H os--Rényi random graph \(G(n,p)\) with \(0 < p < 1\) fixed, \(\kappa\) is equal to the minimum degree \(\delta\) with probability tending to one. On the other hand, Bhattacharya and Mukherjee \cite{bhattacharya_degree_2017} showed that the distribution of \(\delta/\sqrt n\) over inversion graphs approaches the Rayleigh distribution. If it is true that \(\mathbb P(\kappa = \delta) \to 1\) in the class of inversion graphs, then \(\kappa/\sqrt n\) also approaches the Rayleigh distribution for inversion graphs. This would explain the apparent shape of the rows.

\section*{Acknowledgements}

The author is supported by the Swedish Research Council (VR) grant 2022-03875.

\paragraph{Statement on AI use}\addcontentsline{toc}{section}{Statement on AI use}

Large parts of this work were developed in close collaboration with several OpenAI models, namely, GPT 5.6 through 6.1 Sol and GPT 6 Astra. I will outline the contribution of both parties.

I initiated the project when I realized that inversion graphs provide a nice perspective on almost decomposability from \cite{linusson_enumerating_2025}. I proved the \(k = 2\) case of Theorem~\ref{thm:k connected}, and observed that it should generalize to \(\ell\)-chordal graphs. I identified that the notion of an induced ear decomposition would prove it. AI helped find the reference \cite{kriesell_contractible_2000}, and this led immediately to Proposition~\ref{prop:l chordal bound}.

I then used induced ear decompositions to prove (without AI intervention) the following, weaker version of Theorem~\ref{thm:2 connected av}.

\begin{theorem*}
  If every component of \(G_{\sigma}\) is a \(4\)-cycle, then there exists a constant \(C\) depending only on \(\sigma\) such that every \(\pi \in \Av(\sigma)\) whose inversion graph is \(2\)-connected satisfies
  \begin{equation*}
    \inv(\pi) \geq 2|\pi| - C.
  \end{equation*}
\end{theorem*}

My proof of this result was almost identical to the proof of the stronger Theorem~\ref{thm:2 connected av}, with the key difference that it did not use \emph{greedy} ear decompositions.

Next, I asked AI if these two preliminary results generalize to \(k\)-connected inversion graphs. Eventually, it provided me with statements and proofs for both; the former result is Theorem~\ref{thm:k connected}, and the latter is the \(k\)-connected analogue of the statement above. These proofs were very technical and I have not read them in detail, though they were verified in Lean.

I then embarked on a quest to find human-readable proofs of these results, and eventually asked AI for a hint. It suggested using Menger's theorem as in Construction~\ref{const:menger}. Based on this, I was then able to derive proofs for both of the \(k\)-generalizations with the strategy presented in this article.

Finally, I asked if the pattern-avoidance statement could be strengthened to a characterization. Based on computations, AI conjectured Theorem~\ref{thm:2 connected av}, i.e.\ the classification \(\sigma \in \Av(321, 2341, 4123)\). I saw that the claim would follow if one could prove Theorem~\ref{thm:direct sum closure} (closure under direct sums). I asked for hints on how to prove this, and eventually got the suggestion to use greedy ear decompositions and apply the strategy from before. Using this idea, I was able to prove Theorem~\ref{thm:direct sum closure}, and subsequently Theorem~\ref{thm:k connected av}.

The classification of patterns \(\sigma\) in Theorem~\ref{thm:weakly 2-connected} was also conjectured by AI. Simultaneously, it proved the (mostly computational) Lemma~\ref{lem:weakly 2-forcing class}, giving one of the implications in the theorem. The proof written in this paper uses my own version of the argument. I proved the other implication (all of Section~\ref{sec:weak 2-connectivity}) myself.

In summary, the proofs in the final manuscript are mine, but they rely on a few important AI-generated hints.

During the writing of this manuscript, I used AI (both ChatGPT and Claude) for help with programming the figures, finding relevant references, and proofreading. I wrote all of the text, and take full responsibility for the content.

\printbibliography

\clearpage
\appendix
\section{Some constructions} \label{sec:appendix}

This appendix collects some constructions needed for the less interesting implications in Theorems~\ref{thm:k connected av} and \ref{thm:weakly 2-connected}. The first lemma proves that all \(k\)-forcing patterns must belong to \(\Av(321, 2341, 4123)\).

\begin{figure}[t]
  \centering
  \begin{tikzpicture}[scale=0.7]
    \def\ladderedges{
      \foreach \i/\j/\c in {1/2/RoyalBlue, 2/3/BrickRed, 3/4/RoyalBlue, 4/5/BrickRed, 5/6/RoyalBlue}
        \fill[\c!30] (a\i) -- (b\i) -- (b\j) -- (a\j) -- cycle;
      \draw[rpath]
        (a1) -- (a2) -- (a3) -- (a4) -- (a5) -- (a6)
        (b1) -- (b2) -- (b3) -- (b4) -- (b5) -- (b6);
      \foreach \i in {1, ..., 6}
        \draw[rpath] (a\i) -- (b\i);
      \foreach \i in {1, ..., 6}
        \node[dot] at (a\i) {} node[dot] at (b\i) {};
    }

    \foreach \x\y\l in {1/3/a1, 2/5/b2, 3/1/b1, 4/7/a3, 5/2/a2, 6/9/b4,
        7/4/b3, 8/11/a5, 9/6/a4, 10/12/b6, 11/8/b5, 12/10/a6}
      \coordinate (\l) at (\x,\y);
    \ladderedges

    \begin{scope}[shift={(14.45,0)}]
      \foreach \i in {1, ..., 6} {
        \coordinate (a\i) at (0,2.2*\i-1.2);
        \coordinate (b\i) at (2,2.2*\i-1.2);
      }
      \ladderedges
    \end{scope}
  \end{tikzpicture}
  \caption{A diagram of \(\pi_6\) from the proof of Lemma~\ref{lem:not k-forcing} (left); and the structure of the graph (right).}
  \label{fig:pi big}
\end{figure}
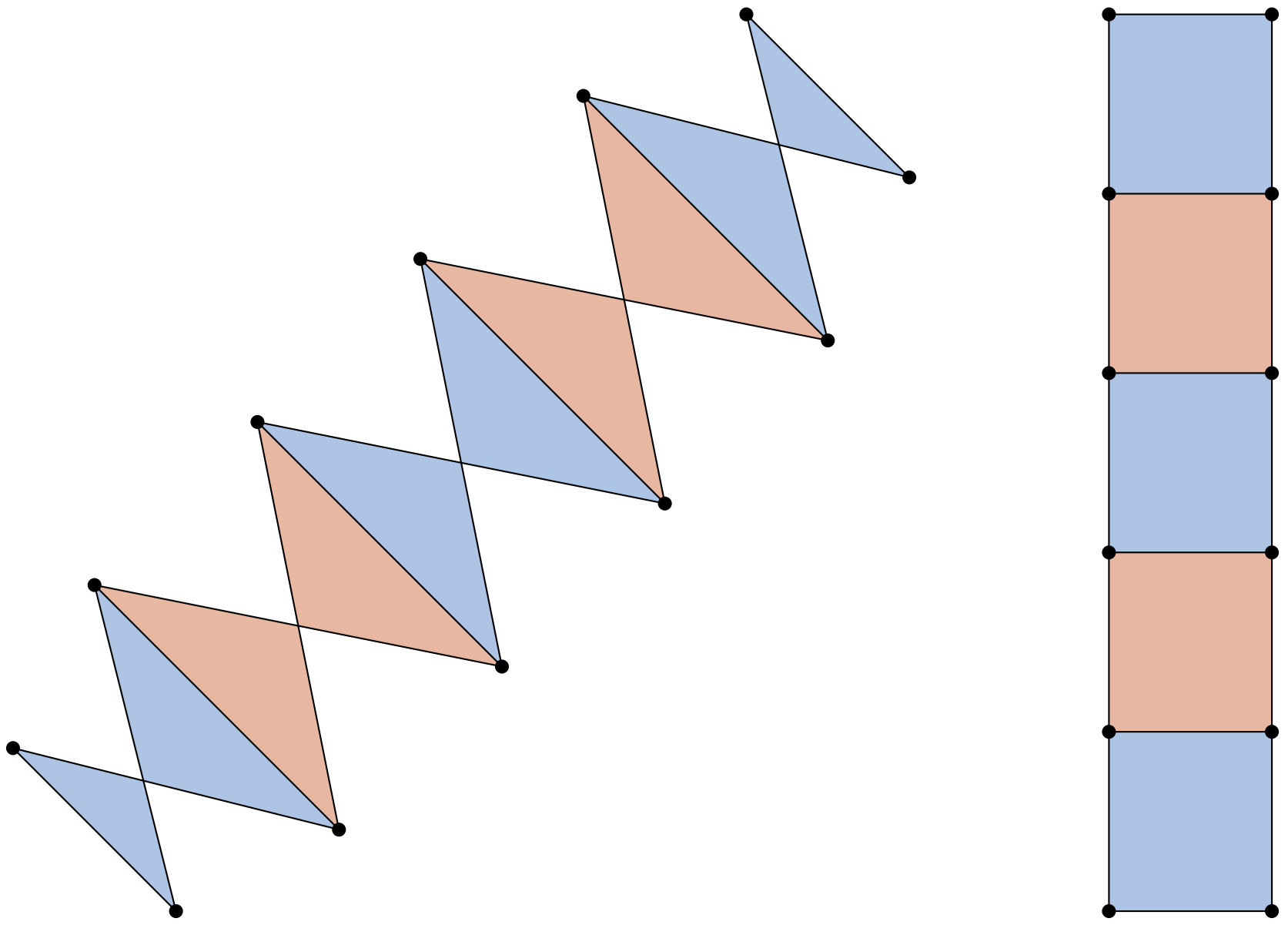

\begin{lemma} \label{lem:not k-forcing}
  The patterns \(321\), \(2341\) and \(4123\) are not \(k\)-forcing for any \(k \geq 2\).
\end{lemma}

For convenience, we define the \emph{\(k\)-deficit}
\begin{equation*}
  \deficit_k(G) = k|V(G)| - |E(G)|.
\end{equation*}

\begin{proof}
  For each \(\sigma\) in the statement, it suffices to find an infinite sequence of permutations \((\pi_n)_n\) in \(\Av(\sigma)\), each with a \(k\)-connected inversion graph, such that 
  \begin{equation*}
    \lim_{n \to \infty} \deficit_k(G_{\pi_n}) = \infty.
  \end{equation*}
  We will start by doing this for \(k = 2\). Since \(4123\) is the inverse of \(2341\), we only need to consider one of the two.

  First, for \(\sigma = 321\), let \(\pi_n\) be the unique permutation of size \(2n\) whose inversion graph is the \emph{ladder graph} \(L_n\) (see Figure~\ref{fig:pi big}); equivalently, \(G_{\pi_n}\) has an induced ear decomposition \(H_0, \ldots, H_{n-2}\) where \(H_0\) is a \(4\)-cycle and every step \(H_i \to H_{i+1}\) adds two new vertices. It is evident that the resulting graph contains no \(3\)-cycles, and it satisfies
  \begin{equation*}
    \deficit_2(G_{\pi_n}) = \deficit_2(H_0) + n - 2 = n + 2 \longrightarrow \infty.
  \end{equation*}

  For \(\sigma = 2341\), slightly more subtlety is called for. We construct permutations \(\tau_n\), with \(\tau_0 = 3412\), inductively as follows (see Figure~\ref{fig:tau big}):
  \begin{itemize}
    \item If \(n \equiv 1\! \pmod 3\), define \(\tau_n\) by inserting a point into \(\tau_{n-1}\), so that the new point is second from the right and second from the top.
    \item If \(n \equiv 2\! \pmod 3\), define \(\tau_n\) by inserting a point into \(\tau_{n-1}\), so that the new point is rightmost and third from the top.
    \item If \(n \equiv 0\! \pmod 3\), define \(\tau_n\) by inserting two new points into \(\tau_{n-1}\) according to Lemma~\ref{lem:inv graph 2 ears} \ref{lem:inv graph 2 ears:top right}.
  \end{itemize}
  Due to the repeating structure of the graph, an inspection of Figure~\ref{fig:tau big} (left) is enough to prove that each \(\tau_n\) avoids \(2341\). For \(n = 3m\) we have
  \begin{equation*}
    \deficit_2(G_{\tau_n}) = \deficit_2(G_{\tau_0}) + m = 4 + m \longrightarrow \infty,
  \end{equation*}
  concluding the proof of the case \(k = 2\).

  \begin{figure}[t]
    \centering
    \begin{tikzpicture}[scale=0.43]
      \def\tauedges{
        \fill[RoyalBlue!30] (a) -- (c) -- (b) -- (p2) -- (p1) -- (d) -- cycle;
        \fill[BrickRed!30] (b) -- (p2) -- (u3) -- (p5) -- (p4) -- (v3) -- cycle;
        \fill[RoyalBlue!30] (u3) -- (p5) -- (u6) -- (p8) -- (p7) -- (v6) -- cycle;
        \fill[BrickRed!30] (u6) -- (p8) -- (u9) -- (p11) -- (p10) -- (v9) -- cycle;
        \draw[rpath]
          (a) -- (c) -- (b) -- (d) -- cycle
          (b) -- (p1) -- (d)
          (b) -- (p2) -- (p1)
          (p2) -- (u3) -- (v3) -- (b)
          (u3) -- (p4) -- (v3)
          (u3) -- (p5) -- (p4)
          (p5) -- (u6) -- (v6) -- (u3)
          (u6) -- (p7) -- (v6)
          (u6) -- (p8) -- (p7)
          (p8) -- (u9) -- (v9) -- (u6)
          (u9) -- (p10) -- (v9)
          (u9) -- (p11) -- (p10);
        \foreach \l in {a, b, c, p1, d, u3, p2, p4, v3, u6, p5, p7, v6, u9, p8, p10, v9, p11}
          \node[dot] at (\l) {};
      }

      \foreach \x\y\l in {1/3/a, 2/7/b, 3/1/c, 4/5/p1, 5/2/d, 6/11/u3, 7/4/p2, 8/9/p4, 9/6/v3, 10/15/u6,
          11/8/p5, 12/13/p7, 13/10/v6, 14/19/u9, 15/12/p8, 16/17/p10, 17/14/v9, 19/16/p11}
        \coordinate (\l) at (\x,\y);
      \tauedges

      \begin{scope}[shift={(23,0)}]
        \foreach \l [count=\i] in {c, b, v3, p4, p5, u6, v9, p10, p11}
          \coordinate (\l) at (0,1+2*\i);
        \foreach \l [count=\i from 0] in {a, d, p1, p2, u3, v6, p7, p8, u9}
          \coordinate (\l) at (3.5,1+2*\i);
        \tauedges
      \end{scope}
    \end{tikzpicture}
    \caption{A diagram of \(\tau_{11}\) from the proof of Lemma~\ref{lem:not k-forcing} (left); and the structure of the graph (right).}
    \label{fig:tau big}
  \end{figure}
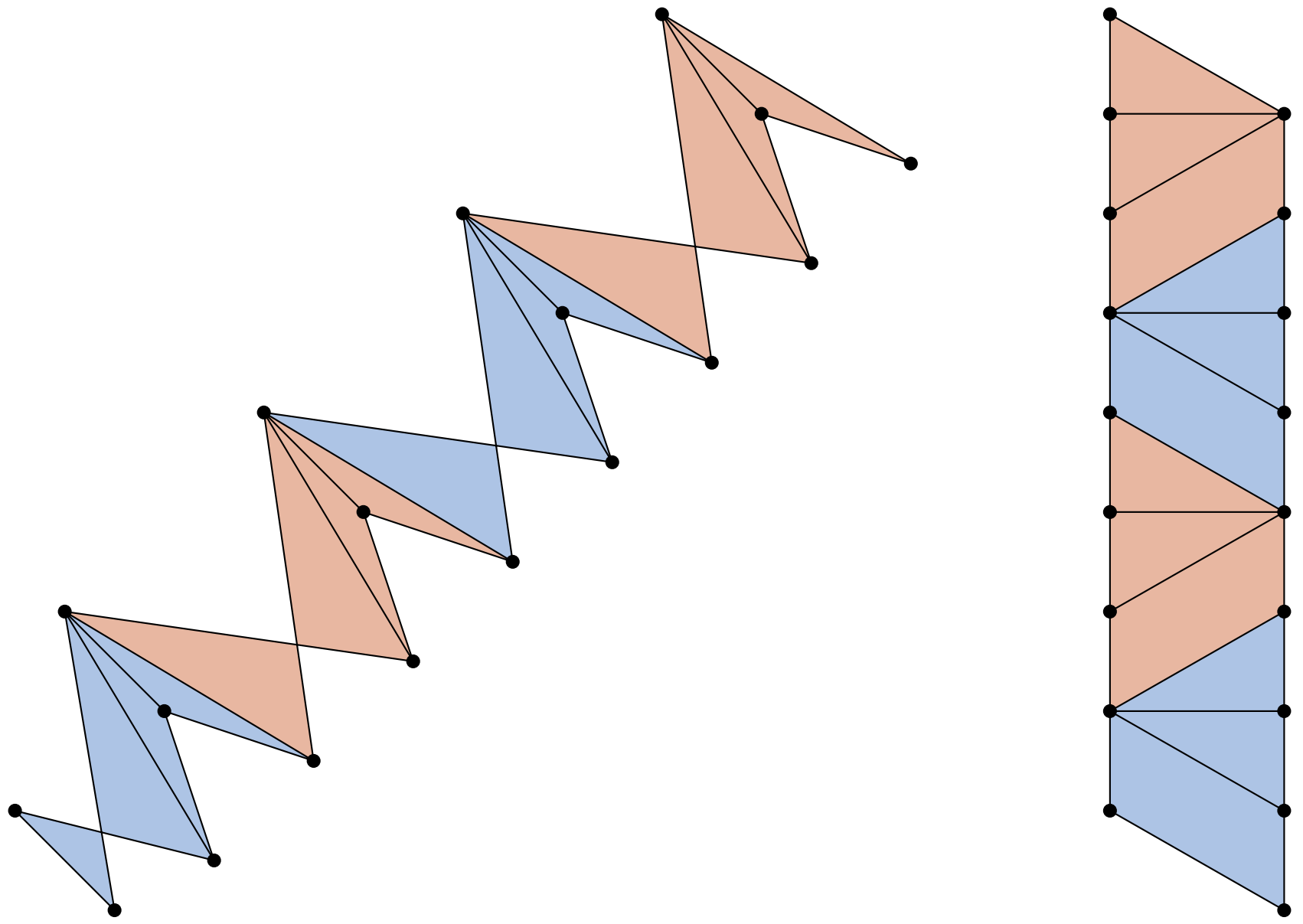

  Fix now an arbitrary \(k \geq 2\); we will generalize our previous constructions. Denote by \(L^k_n\) the graph on \(kn\) vertices obtained by inflating each vertex in one pillar of the ladder graph \(L_n\) into an independent set of size \(k-1\). (When a vertex \(v\) is inflated into \(\{v_1, \ldots, v_{k-1}\}\), each \(v_i\) gets the same set of neighbors as the old vertex \(v\).) Figure~\ref{fig:pi prime big} (right) shows an example for \(k = 3\). There exists a permutation \(\pi^k_n\) whose inversion graph is isomorphic to \(L^k_n\), as Figure~\ref{fig:pi prime big} (left) shows.
  
  Since \(L^k_n\) is bipartite, \(\pi^k_n\) must avoid \(321\). It is also clear that \(L^k_n\) is \(k\)-connected. Furthermore,
  \begin{equation*}
        |E(L^k_n)| = |E(L^k_{n-1})| + (k-1)^2 + k = (k^2 - k + 1)n - (k^2 - 2k + 2),
  \end{equation*}
  and thus
  \begin{equation*}
    \deficit_k(L^k_n) = k^2n - |E(L^k_n)| = (k - 1)n + (k - 1)^2 + 1 \longrightarrow \infty.
  \end{equation*}

  \begin{figure}[t]
    \centering
    \begin{tikzpicture}[scale=0.67]
      \def\inflatedladderedges{
        \foreach \i in {1, ..., 6}
          \draw[black!20, line width=12pt, line cap=round] (a\i1) -- (a\i2);
        \draw[rpath]
          (b1) -- (b2) -- (b3) -- (b4) -- (b5) -- (b6);
        \foreach \i/\j in {1/2, 2/3, 3/4, 4/5, 5/6}
          \draw[rpath] (a\i1) -- (a\j1) -- (a\i2) -- (a\j2) -- cycle;
        \foreach \i in {1, ..., 6}
          \draw[rpath] (a\i1) -- (b\i) -- (a\i2);
        \foreach \i in {1, ..., 6}
          \node[dot] at (a\i1) {} node[dot] at (a\i2) {} node[dot] at (b\i) {};
      }

      \foreach \x\y\l in {3/1/b1, 2/5/b2, 7/4/b3, 6/9/b4, 11/8/b5, 10/12/b6}
        \coordinate (\l) at (\x,\y);
      \foreach \x\y\l in {1/3/a1, 5/2/a2, 4/7/a3, 9/6/a4, 8/11/a5, 12/10/a6} {
        \coordinate (\l1) at (\x-0.25,\y-0.25);
        \coordinate (\l2) at (\x+0.25,\y+0.25);
      }
      \inflatedladderedges

      \begin{scope}[shift={(14.7,0)}]
        \foreach \i in {1, ..., 6} {
          \coordinate (a\i1) at (0,2.1*\i-0.85);
          \coordinate (a\i2) at (0.8,2.1*\i-0.6);
          \coordinate (b\i) at (2.8,2.1*\i-1.1);
        }
        \inflatedladderedges
      \end{scope}
    \end{tikzpicture}
    \caption{A diagram of \(\pi^k_6\) for \(k = 3\) from the proof of Lemma~\ref{lem:not k-forcing} (left); and the structure of the graph (right).}
    \label{fig:pi prime big}
  \end{figure}
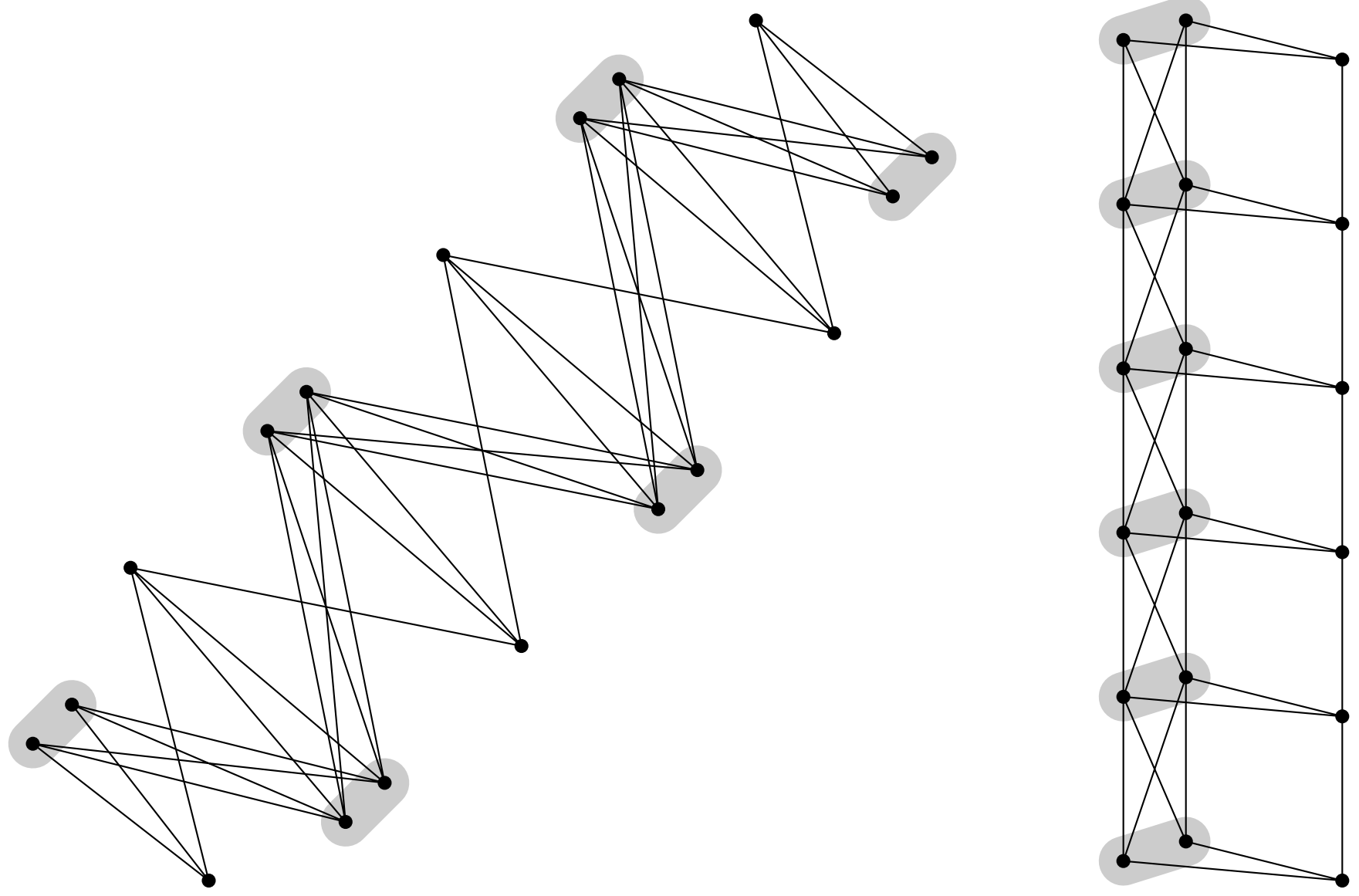
  
  The case \(\sigma = 2341\) remains. We define
  \begin{equation*}
    \tau^k_n = \delta_{k-2} \ominus \tau_n,
  \end{equation*}
  where \(\delta_{\ell}\) denotes the decreasing permutation \(\ell \ldots 2 1\), and \(\tau_n\) is as before.

  The inversion graph of \(\tau^k_n\) is evidently \(k\)-connected, and it avoids \(2341\), since an occurrence of \(2341\) cannot use any of the new \(k-2\) points. Also, recall that if \(n = 3m\), then \(\deficit_2(G_{\tau_n}) = 4 + m\). Using this, a short calculation gives
  \begin{equation*}
    \deficit_k(G_{\tau_n^k}) = 4 + m + \frac{(k-2)(k+3)}{2} \longrightarrow \infty.
  \end{equation*}
  This concludes the proof.
\end{proof}

\begin{remark} \label{rmk:weaker forcing example}
  If we write \(m = |\tau_n|\) with \(\tau_n\) as in the proof above, our observation about the deficit of \(G_{\tau_n}\) can be rewritten as
  \begin{equation*}
    \inv(\tau_n) = \left\lceil\frac{7}{4} m\right\rceil - 3.
  \end{equation*}
  Therefore, this family of permutations provides an example of a linear coefficient between the usual \(\frac{3}{2}\) for \(2\)-connected inversion graphs and \(2\) from Theorem~\ref{thm:k connected av}.
\end{remark}

The next result is used to prove the implication \(\ref{thm:weakly 2-connected:forcing} \implies \ref{thm:weakly 2-connected:patterns}\) in Theorem~\ref{thm:weakly 2-connected}. Concretely, we want to show that every pattern not contained in \(\Sub(\Gamma B)\) from the theorem statement fails to be strongly \(2\)-forcing. It suffices to find a collection \(F\) such that \(\Sub(\Gamma B) = \Av(F)\), and prove that no pattern in \(F\) is strongly \(2\)-forcing.

\begin{figure}[t]
  \centering
  \begin{tikzpicture}[scale=1.1]
    \coordinate (sw) at (0,0);
    \coordinate (ne) at (6,6);
    \coordinate (a) at (2.2,3.8);
    \coordinate (b) at (3.8,2.2);

    \fill[black!20]
      (sw) rectangle (a |- b)
      (b |- a) rectangle (ne);
    \draw (sw) rectangle (ne);
    \draw[rpath]
      (a |- sw) -- (a |- ne)
      (b |- sw) -- (b |- ne)
      (a -| sw) -- (a -| ne)
      (b -| sw) -- (b -| ne);

    \foreach \l in {a, b}
      \node[dot] at (\l) {};
    \node[above left] at (a) {\(\pi(i)\)};
    \node[below right] at (b) {\(\pi(j)\)};
    \node at (1.1,1.1) {\(\mathrm{SW}\)};
    \node at (4.9,4.9) {\(\mathrm{NE}\)};
  \end{tikzpicture}
  \caption{The diagram of \(G_{\pi}\) in the proof of Lemma~\ref{lem:weakly 2-forcing class}.}
  \label{fig:weakly 2-forcing class lemma}
\end{figure}
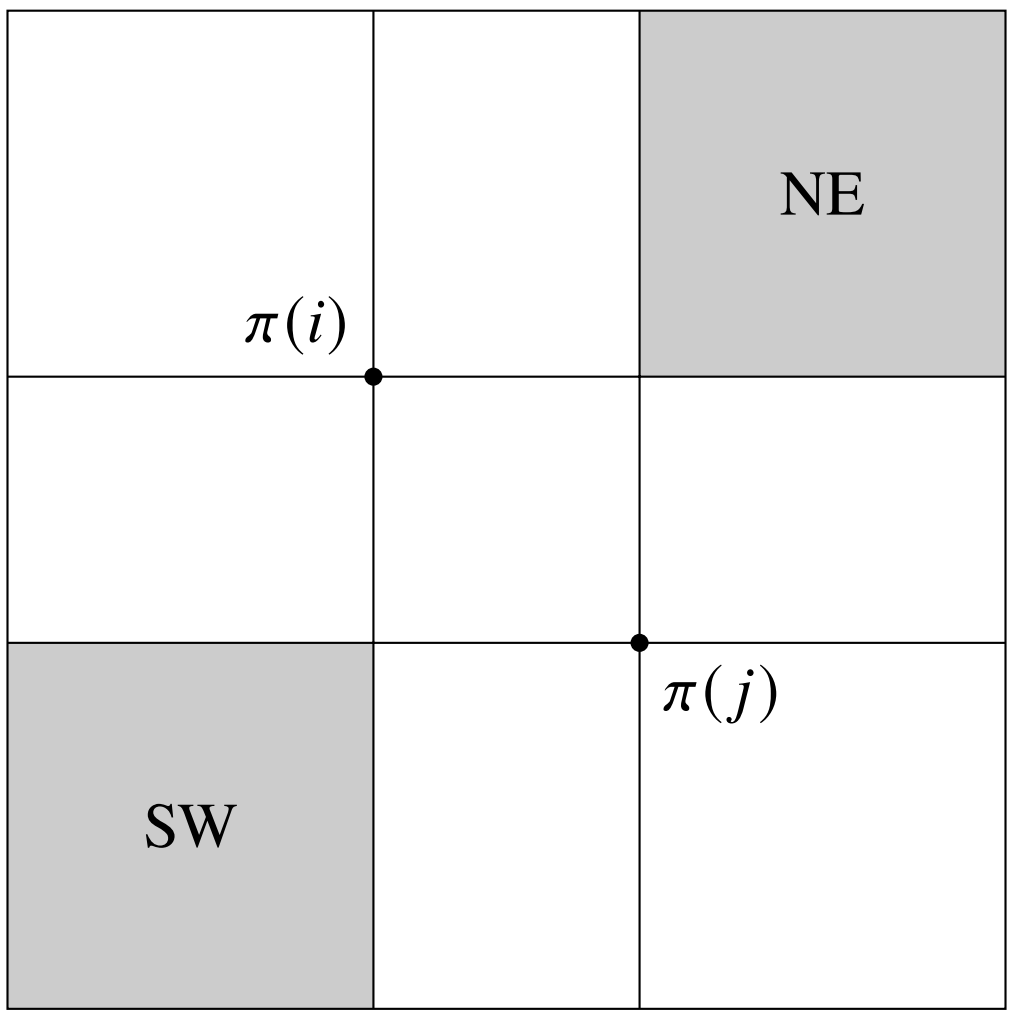

\begin{lemma} \label{lem:weakly 2-forcing class}
  Let \(B\) be as in Theorem~\ref{thm:weakly 2-connected} \ref{thm:weakly 2-connected:patterns}, and let
  \begin{equation*}
    F = \left\{
    \begin{aligned}
      &321, \\
      &2341, 4123, 3412, \\
      &12435, 12453, 12534, 13245, 13254, \\
      &13425, 13524, 14235, 14253, 21435, \\
      &21453, 21534, 23145, 23154, 24135, \\
      &24153, 31245, 31254, 31425, 31524 
    \end{aligned}
    \right\}.
  \end{equation*}
  Then \(\Sub(\Gamma B) = \Av(F)\). Furthermore, no pattern in \(F\) is strongly \(2\)-forcing.
\end{lemma}

\begin{remark}
  \(F\) contains \(24\) patterns. The first three, i.e.\ \(321\), \(2341\) and \(4123\), are familiar from Theorem~\ref{thm:k connected av}. The last \(21\) are `noise' introduced by the more contrived notion of weak \(2\)-connectivity.
\end{remark}

\begin{proof}[Proof of Lemma~\ref{lem:weakly 2-forcing class}]
  We will begin by showing that \(\Av(F) \subseteq \Sub(\Gamma B)\). For \(n \leq 5\), it is easy to compute \(\Av_n(F)\) and verify that this holds (with equality).
  
  Therefore, it suffices to prove that for \(n \geq 6\), the class \(\Av_n(F)\) contains only the four patterns \(\id_n\), \({\id_{n-2}} \oplus 21\), \(21 \oplus {\id_{n-2}}\) and \(21 \oplus {\id_{n-4}} \oplus 21\); clearly, these are the only patterns of size \(n\) in \(\Sub(\Gamma B)\).

  To see this, suppose that \(\Av_n(F)\) contains some other permutation \(\pi\). This permutation must have an inversion \(i < j\), \(\pi(i) > \pi(j)\), that does \emph{not} satisfy
  \begin{equation*}
    \{i, j\} = \{\pi(i), \pi(j)\} \in \{\{1, 2\},  \{n-1, n\}\}.
  \end{equation*}
  We will show that there exists a pattern \(\tau\) of size \(5\) in \(\pi\) such that
  \begin{equation*}
    \tau \notin \Av_5(F) = \{12345, 12354, 21345, 21354\}.
  \end{equation*}
  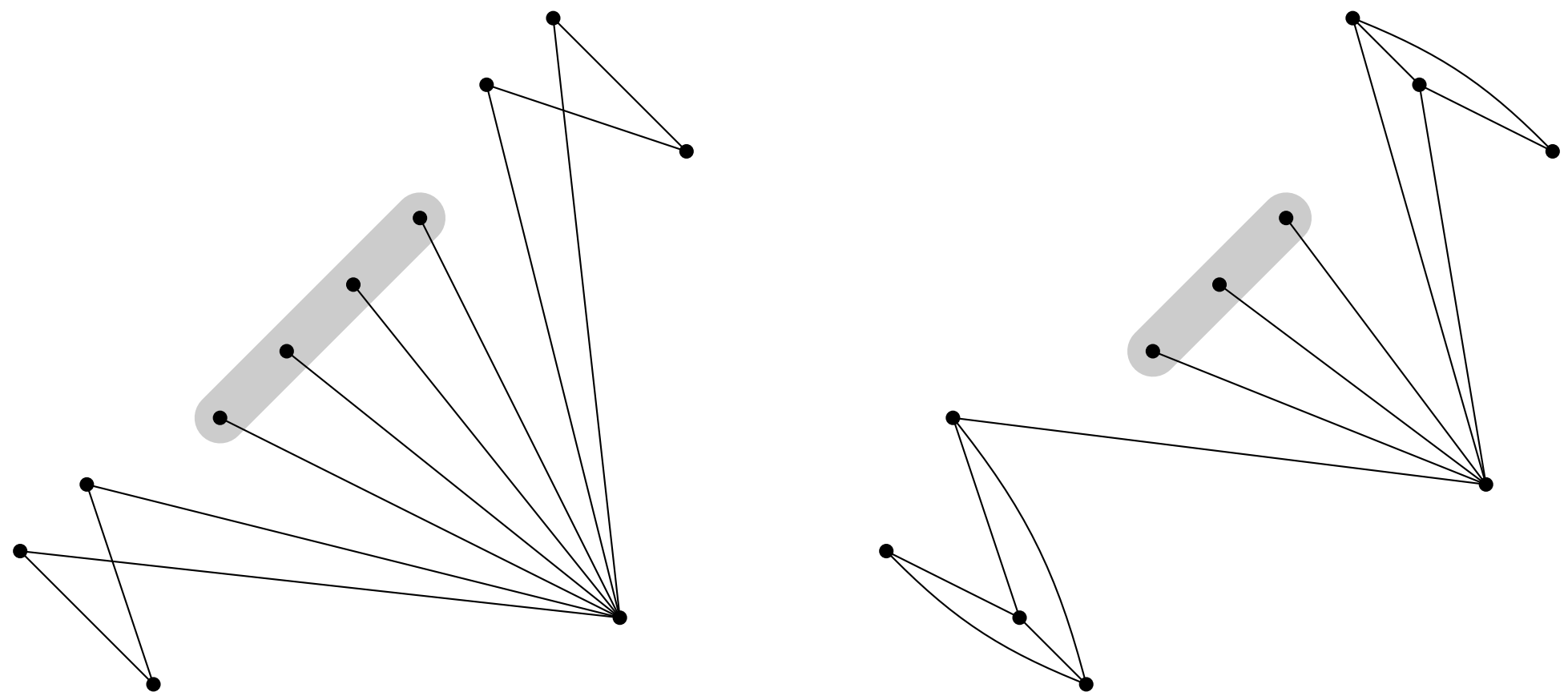
\begin{figure}[t]
    \centering
    \begin{tikzpicture}[scale=0.55]
      \foreach \v [count=\i] in {3, 4, 1, 5, 6, 7, 8, 10, 11, 2, 9}
        \coordinate (\v) at (\i,\v);
      \draw[black!20, line width=12pt, line cap=round] (5) -- (8);
      \foreach \a/\b in {3/1, 3/2, 4/1, 4/2, 5/2, 6/2, 7/2, 8/2, 10/2, 10/9, 11/2, 11/9}
        \draw[rpath] (\a) -- (\b);
      \foreach \v in {1, ..., 11}
        \node[dot] at (\v) {};
      \node at (6,-0.3) {\(\alpha\)};

      \begin{scope}[shift={(13,0)}]
        \foreach \v [count=\i] in {3, 5, 2, 1, 6, 7, 8, 11, 10, 4, 9}
          \coordinate (\v) at (\i,\v);
        \draw[black!20, line width=12pt, line cap=round] (6) -- (8);
        \foreach \a/\b in {3/2, 5/2, 5/4, 2/1, 6/4, 7/4, 8/4, 11/10, 11/4, 10/4, 10/9}
          \draw[rpath] (\a) -- (\b);
        \draw[rpath]
          (3) to[bend right=12] (1)
          (5) to[bend left=12] (1)
          (11) to[bend left=12] (9);
        \foreach \v in {1, ..., 11}
          \node[dot] at (\v) {};
        \node at (6,-0.3) {\(\beta\)};
      \end{scope}
    \end{tikzpicture}
    \caption{Inversion graphs of the permutations \(\alpha\) and \(\beta\) for \(n = 11\).}
    \label{fig:alpha beta}
  \end{figure}

  Consider the nine regions in the diagram of \(G_{\pi}\) cut out by the horizontal and vertical lines through \(\pi(i)\) and \(\pi(j)\), visualized in Figure~\ref{fig:weakly 2-forcing class lemma}.

  By our earlier observation regarding the inversion \((i, j)\), it is not possible that all remaining points of \(G_{\pi}\) are in the region SW (southwest of both \(\pi(i)\) and \(\pi(j)\)). Similarly, they cannot all be in NE. We have two cases:
  \begin{itemize}
    \item If any of the white regions contains a point \(v\), then let \(\tau\) be a length-five pattern of \(\pi\) containing all of \(v\), \(\pi(i)\) and \(\pi(j)\).

    \item Otherwise, SW and NE must both be nonempty. Let \(\tau\) be a length-five pattern containing a point from each, and both \(\pi(i)\) and \(\pi(j)\).
  \end{itemize}
  In the former case, \(\tau\) must contain \(231\), \(312\) or \(321\). In the latter, it must contain \(1324\). All elements of \(\Av_5(F)\) avoid these patterns, which is a contradiction: \(\Av(F)\) is closed under pattern containment. We conclude that \(\Av(F) \subseteq \Sub(\Gamma B)\).

  For the reverse inclusion, it suffices to examine the four permutations in \(\Sub(\Gamma B)\) of size \(n \geq 5\), and observe that they avoid all the patterns in \(F\). Thus, \(\Av(F) = \Sub(\Gamma B)\).
  
  The second part of the claim was that no pattern in \(F\) is strongly \(2\)-forcing. To do this, we shall consider two special permutations:
  \begin{align*}
    \alpha &= [3,\ 4,\ 1,\ \ 5,\ 6,\ \ldots,\ n-3,\ \ n-1,\ n,\ 2,\ n-2], \\
    \beta &= [3,\ 5,\ 2,\ 1,\ \ 6,\ 7,\ \ldots,\ n-3,\ \ n,\ n-1,\ 4,\ n-2].
  \end{align*}
  Their inversion graphs are shown in Figure~\ref{fig:alpha beta}, with the expandable increasing sequences highlighted.

  It is clear that both are weakly \(2\)-connected, and they have \(n+1\) and \(n+3\) inversions, respectively. The proof is complete if, up to symmetry, each element of \(F\) is avoided by one of the permutations \(\alpha\) and \(\beta\) for every \(n \geq 8\). This is easy to check by hand.
\end{proof}

\end{document}